\documentclass[11pt,oneside,leqno]{amsart}
\usepackage{amsxtra}
\usepackage{amsopn}
\usepackage{color}
\usepackage{amsmath,amsthm,amssymb}
\usepackage{amscd}
\usepackage{amsfonts}
\usepackage{latexsym}
\usepackage{verbatim}

\definecolor{verde}{rgb}{0.5,.7,.2}

\theoremstyle{plain}
\newtheorem{theorem}{Theorem}[section]
\newtheorem{definition}[theorem]{Definition}
\newtheorem{lemma}[theorem]{Lemma}
\newtheorem{proposition}[theorem]{Proposition}
\newtheorem{corollary}[theorem]{Corollary}
\newtheorem{remark}[theorem]{Remark}
\newtheorem{example}[theorem]{Example}
\newtheorem{question}[theorem]{QUESTION}
\newtheorem{remark-question}[section]{Remark-Question}
\newtheorem{conjecture}[section]{Conjecture}

\newcommand\fra{{\mathfrak a}} 

\newcommand\frg{{\mathfrak g}}
\newcommand\frh{{\mathfrak h}}

\newcommand*{\longlongrightarrow}{\ensuremath{\joinrel\relbar\joinrel\relbar\joinrel\relbar\joinrel\relbar\joinrel\rightarrow}}

\begin{document}
\title[Construction of Lie algebras with special $\mathrm{G}_2$-structures]{Construction of Lie algebras with special $\mathrm{G}_2$-structures}
\date{\today}
\author{V\'ictor Manero}
\maketitle


\begin{abstract}
We give a method to obtain new 7-dimensional Lie algebras endowed with closed and coclosed $\mathrm{G}_2$-structures starting from 6-dimensional Lie algebras with symplectic half-flat $\mathrm{SU}(3)$-structures and half-flat $\mathrm{SU}(3)$-structures, respectively
. Finally, we describe all the 7-dimensional Lie algebras with a closed $\mathrm{G}_2$-structure that are obtained with this method from the 6-dimensional solvable Lie algebras admitting a symplectic half-flat $\mathrm{SU}(3)$-structure.
 \end{abstract}


\begin{section}*{Introduction}
An $\mathrm{SU}(n)$-structure on a Lie algebra $\mathfrak{h}$ of dimension $2n$, consists in a triple $(g, J, \Psi)$ such that $(g, J)$ is an almost Hermitian structure on $\mathfrak{h}$, and $\Psi=\psi_++ i \, \psi_-$ is a complex volume $(n,0)$-form, satisfying 
\begin{equation*}
(-1)^{n(n-1)/2}\,\Big(\frac{i}{2}\Big)^{n}\,\Psi\wedge \overline{\Psi}=\frac{1}{n!} \, \omega^n,
\end{equation*}
with $\overline{\Psi}$ the complex form obtained by conjugation of $\Psi$, and $\omega$ the K\"ahler form associated to $(g, J)$. In what follows we will consider $\mathrm{SU}(3)$-structures on 6-dimensional Lie algebras. 

The existence of an $\mathrm{SU}(3)$-structure on a Lie algebra $\mathfrak{h}$ can also be described by the presence of a pair of forms, namely, $(\omega, \psi_+)\in \Lambda^2\mathfrak{h}^* \times \Lambda^3\mathfrak{h}^*$ such that describe a metric as
 \begin{equation*}
 g(X,Y)\,\omega^3=-3\,\iota_X\omega\wedge\iota_Y(\psi_+)\wedge\psi_+,
 \end{equation*}
 with $X, Y \in \mathfrak{h}$ and $\iota_X$ denoting the contraction by $X$. We can also recover its compatible almost complex structure as it is described in \cite{CLSS} 
  \begin{equation*}
 (J_{\psi_+}^*\alpha)(X)\,\omega^3=\alpha\wedge \iota_X\psi_+\wedge \psi_+,
 \end{equation*}
 or, equivalently,
 \begin{equation*}
 \alpha(JX)=-J^*\alpha(X),
 \end{equation*}
for any 1-form $\alpha$ on $\mathfrak{h}^*$.

Also, if $(g, J, \Psi)$ is an $\mathrm{SU}(3)$-structure on a Lie algebra $\mathfrak{h}$ we may choose an orthonormal frame $\{e_1,\dots, e_6\}$ such that the almost complex structure $J$ is $J^*e^1=e^2$, $J^*e^3=e^4$ and $J^*e^5=e^6$ with $\{e^1,\dots, e^6\}$ an orthonormal basis dual to $\{e_1,\dots, e_6\}$. Therefore, the K\"ahler form $\omega$ and the complex volume form $\Psi$ can be written as
\begin{equation}\label{SU(3)}
 \omega=e^{12}+e^{34}+e^{56}, \qquad \Psi=(e^1+ i \, e^2)\wedge (e^3+ i \, e^4)\wedge (e^5+ i \, e^6),
 \end{equation}
 where, with the usual notation of the related literature, we write 
  $e^{ij}$ for the wedge product $e^{i}\,\wedge\,e^{j}$,   $e^{ijk}=e^{i}\,\wedge\,e^{j}\,\wedge\,e^{k}$, and so on.
  Thus, 
 \begin{equation*}
 \psi_+=e^{135}-e^{146}-e^{236}-e^{245}, \qquad \text{ and } \qquad \psi_-=-e^{246}+e^{235}+e^{145}+e^{136}.
 \end{equation*}
 In \cite{GH}, Gray and Hervella prove that there exist sixteen different classes of almost Hermitian structures attending to the behavior of the covariant derivative of its K\"ahler form. Equivalently, the different classes of $\mathrm{SU}(n)$-structures can be defined in terms of the forms $\omega, \psi_+$ and $\psi_-$. In particular we are interested in two classes of $\mathrm{SU}(3)$-structures which were defined respectively in \cite{CS} and \cite{TV} as follows: 

\medskip

\begin{itemize}
\item $(g, J, \Psi)$ is a \emph{half-flat} $\mathrm{SU}(3)$-structure iff $d\omega^{2}=d\psi_+=0$;
\medskip
\item $(g, J, \Psi)$ is a \emph{symplectic half-flat} $\mathrm{SU}(3)$-structure iff $d\omega=d\psi_+=0$.
\end{itemize}

\medskip

A classification of half-flat $\mathrm{SU}(3)$-structures on nilpotent Lie algebras is done in \cite{Co}.  In \cite{FSH1} a similar work for indecomposable solvable Lie algebras has been stablished. The existence of symplectic half-flat $\mathrm{SU}(3)$-structures on nilpotent Lie algebras is studied in \cite{CT} and the complete study of these structures on solvable Lie algebras is obtained in \cite{FMOU}. 

\bigskip

A $\mathrm{G}_2$-structure on a 7-dimensional Lie algebra $\mathfrak{g}$ is defined by a 3-form $\varphi$ (called the fundamental form) on $\mathfrak{g}$ such that 
 \begin{equation*}\label{G2-metric}
g_{\varphi}(X,Y)\,vol =\frac{1}{6}\, \iota_X\varphi\wedge\iota_Y\varphi\wedge\varphi,
\end{equation*}
defines a Riemannian metric with $X, Y \in \mathfrak{g}$ and $vol$ denoting the volume form.
With respect to some orthonormal basis of 1-forms $\{e^1,\dots, e^7\}$ on $\mathfrak{g}$ the fundamental form can be written as
\begin{equation}\label{G2form}
\varphi=e^{127}+e^{347}+e^{567}+e^{135}-e^{146}-e^{236}-e^{245}.
\end{equation}
It can also be defined the 4-form $\ast \varphi$, where $\ast$ denotes the Hodge star operator associated to $g_{\varphi}$. Therefore, respect to the basis $\{e^{1}, \dots, e^7\}$ of 1-forms of $\mathfrak{g}$ in which the fundamental form is described by \eqref{G2form} the 4-form can be described as
 \begin{equation*}
\varphi=e^{1234}+e^{1256}+e^{1234}-e^{2467}+e^{2357}+e^{1457}+e^{1367}.
\end{equation*}

In \cite{FG}, Fern\'andez and Gray prove that there exist sixteen different classes of $\mathrm{G}_2$-structures attending to the behavior of the covariant derivative of its  fundamental form. In particular we will be interested in two different classes of $\mathrm{G}_2$-structures which are described as follows:
\medskip

\begin{itemize}
\item $\varphi$ is an \emph{almost parallel} or \emph{closed} $\mathrm{G}_2$-structure iff $d\varphi=0$;
\medskip
\item $\varphi$ is a \emph{semiparallel} or \emph{coclosed} $\mathrm{G}_2$-structure iff $d\ast\varphi=0$.
\end{itemize}

\medskip
A classification of closed $\mathrm{G}_2$-structures on nilpotent Lie algebras has been recently obtained in \cite{CF}.

$\mathrm{SU}(3)$-structures and $\mathrm{G}_2$-structures are closely related. In fact, if $(N^6, \omega, \psi_+)$ is a 6-dimensional manifold endowed with an $\mathrm{SU}(3)$-structure then the 3-form 
\begin{equation}\label{G2formfromSU(3)}
\varphi=\omega\wedge dt+\psi_+,
\end{equation}
defines a $\mathrm{G}_2$-structure on the 7-dimensional manifold $M^7=N^6\times S^1$ where $t$ denotes the coordinate in $S^1$.

Concerning the relation between special $\mathrm{SU}(3)$-structures and special $\mathrm{G}_2$-structures, if the $\mathrm{SU}(3)$-structure $(\omega, \psi_+)$ on $N^6$ is symplectic half-flat clearly the $\mathrm{G}_2$-structure defined by \eqref{G2formfromSU(3)} constitutes a closed $\mathrm{G}_2$-structure on $M^7$. Equivalently, if the $\mathrm{SU}(3)$ manifold $(N^6, \omega, \psi_+)$ is half-flat the 3-form 
\begin{equation}\label{G2coformfromSU(3)}
\varphi=\omega\wedge dt-\psi_-,
\end{equation}
is such that 
\begin{equation*}
\ast\varphi=\frac{1}{2}\omega\wedge \omega+ \psi_+\wedge dt,
\end{equation*}
and therefore defines a coclosed $\mathrm{G}_2$-structure on the 7-dimensional manifold $M^7=N^6\times S^1$ where $t$ is the coordinate on $S^1$.

Regarding the converse, Cabrera in \cite{Ca} shows that if $N^6$ is an orientable hypersurface of a $\mathrm{G}_2$ manifold $(M^7, \varphi)$ then the pair of forms on $N^6$ defined by,
\begin{equation*}
\omega=\iota_U\varphi \qquad and \qquad \psi_+=\pi^*\varphi,
\end{equation*}
with $U$ the unitary vector field of $M^7$ normal to $N^6$ and $\pi$ the projection of $M^7$ onto $N^6$, describe an $\mathrm{SU}(3)$-structure on $N^6$.

\medskip

If we focus our attention on Lie algebras it is clear that the presence of a symplectic half-flat structure namely $(\omega, \psi_+)$ on a 6-dimensional Lie algebra $\mathfrak{h}$, defines a closed $\mathrm{G}_2$-structure on $\mathfrak{g}=\mathfrak{h}\oplus \mathbb{R}$ defined as in \eqref{G2formfromSU(3)}. Equivalently, if the $\mathrm{SU}(3)$-structure $(\omega, \psi_+)$ on $\mathfrak{h}$ is half-flat the $\mathrm{G}_2$-structure defined on $\mathfrak{g}=\mathfrak{h}\oplus \mathbb{R}$ by \eqref{G2coformfromSU(3)} is coclosed. These constructions of 7-dimensional Lie algebras endowed with closed and coclosed $\mathrm{G}_2$-structures as the direct sum of 6-dimensional Lie algebras with symplectic half-flat $\mathrm{SU}(3)$-structures or half-flat ones plus a 1-dimensional abelian Lie algebra are well known. In the present work we generalize this construction. This fact allows to obtain new examples of 7-dimensional Lie algebras endowed with closed and coclosed $\mathrm{G}_2$-structures.
Thus, provided the existence of a lattice we can construct new compact solvmanifolds endowed with special $\mathrm{G}_2$-structures.
\medskip

In Proposition \ref{G2cerrada} we describe how to obtain a 7-dimensional Lie algebra of the form
\begin{equation}\label{algebra}
\frg=\frh\oplus_D\mathbb{R}
\end{equation}  
endowed with a closed $\mathrm{G}_2$-structure from a 6-dimensional Lie algebra $\frh$ with a symplectic half-flat $\mathrm{SU}(3)$-structure, where $D$ denotes a derivation of $\frh$, what constitutes a generalization of some of the results obtained in \cite{Fr}. Then, in Example \ref{ejemplo1} we give a concrete example of this construction. Concerning the converse, in Proposition \ref{G2cerradainversa} we show how a closed $\mathrm{G}_2$-structure on a 7-dimensional Lie algebra of the form \eqref{algebra}, with $D$ a particular type of derivation, describes a symplectic half-flat $\mathrm{SU}(3)$-structure on the 6-dimensional subalgebra $\frh$. An example of this construction is described in Example \ref{ejemplo2}.

Section 2 is devoted to an equivalent study but considering coclosed $\mathrm{G}_2$-structures and half-flat $\mathrm{SU}(3)$-structures. In particular, in Proposition 2.1 we describe how to obtain 7-dimensional Lie algebras endowed with a coclosed $\mathrm{G}_2$-structure from 6-dimensional Lie algebras with half-flat $\mathrm{SU}(3)$-structures. A concrete example of this construction is given in Example \ref{ejemplo3}. Regarding the converse of Proposition \ref{G2cocerrada} we show in Proposition \ref{G2cocerradainversa} how a coclosed $\mathrm{G}_2$-structure on a 7-dimensional Lie algebra $\frg$ of the form \eqref{algebra}, with $D$ a particular type of derivation, describes a half-flat $\mathrm{SU}(3)$-structure on the subalgebra $\frh$. Example \ref{ejemplo4} shows the use of Proposition \ref{G2cocerradainversa}. 
\end{section}


\begin{section}{Lie algebras with a closed $\mathrm{G}_2$-structure}\label{Closed}
We show that if a 6-dimensional symplectic half-flat Lie algebra is endowed with a particular type of derivation, 
then one can construct a Lie algebra with a closed $\mathrm{G}_2$ form. If $\mathfrak{h}$ is a 6-dimensional Lie algebra, and $D$ a derivation of $\mathfrak{h}$, the vector space
\begin{equation*}
\mathfrak{g}=\mathfrak{h}\oplus_{D} \mathbb{R} \xi
\end{equation*}
is a Lie algebra with the Lie bracket given by 
\begin{equation*}
[U,V]=[U,V]|_{\mathfrak{h}}, \qquad [\xi, U]=D(U),
\end{equation*}
for any $U,V \in \mathfrak{h}$.

We recall that a closed $\mathrm{G}_2$ form on a real Lie algebra $\mathfrak{g}$ of dimension 7 is a closed 3-form $\varphi$
on $\mathfrak{g}$ such that can be written as in \eqref{G2form} with respect to some basis 
$\{e^1,\dotsc, e^7\}$ of the dual space of $\mathfrak{g}$. 

Let $(\omega, \psi_{+})$ be a symplectic half-flat structure on $\mathfrak h$. Thus, it defines an almost complex structure $J$, and as it is mention on \cite{BV} 
this allows to obtain a real representation of the complex matrices as
\begin{equation*}
\rho:  \mathfrak{gl}(3, \mathbb{C}) \longlongrightarrow  \mathfrak{gl}(6, \mathbb{R}).
\end{equation*}
Then, if $A\in \mathfrak{gl}(3, \mathbb{C})$, $\rho(A)$ is the matrix $(B_{ij})_{i,j=1}^3$ with 
\begin{equation*}
B_{ij}=\left(
\begin{array}{cc}
 Re A_{ij} & Im A_{ij}   \\
 -Im A_{ij} & Re A_{ij}  \\
\end{array}
\right),
\end{equation*}
where $A_{ij}$ is the $(i,j)$ component of $A$.\\
In particular, the real representation of $\mathfrak{sl}(3,\mathbb{C})$ (complex matrices without trace) is given by
\begin{equation}\label{sl(3,C)}
\mathfrak{sl}(3, \mathbb{C})=\left\{
\left(
\begin{array}{cc|cc|cc}
 a_{1,1} & a_{1,2} & a_{1,3} & a_{1,4} & a_{1,5} & a_{1,6}  \\
 -a_{1,2} & a_{1,1} & -a_{1,4} & a_{1,3} & -a_{1,6} & a_{1,5}  \\\hline
 a_{3,1} & a_{3,2} & a_{3,3} & a_{3,4} & a_{3,5} & a_{3,6}  \\
 -a_{3,2} & a_{3,1} & -a_{3,4} & a_{3,3} & -a_{3,6} & a_{3,5}  \\\hline
 a_{5,1} & a_{5,2} & a_{5,3} & a_{5,4} & -a_{1,1}-a_{3,3} & -a_{1,2}-a_{3,4}  \\
 -a_{5,2} & a_{5,1} & -a_{5,4} & a_{5,3} & a_{1,2}+a_{3,4} & -a_{1,1}-a_{3,3}  \\
\end{array}
\right), \text{ with } a_{i,j}\in \mathbb{R}
\right\}.
\end{equation}
On the other hand, the $\mathrm{SU}(3)$-structure on $\mathfrak{h}$ guarantees the existence of certain basis, namely $\{e^1,\dots, e^6\}$ of $\mathfrak{h}^*$, in which $\omega, \, \psi_+$ and $\psi_-$ have the canonical expression 
\begin{equation*} \label{canonicalexp}
\begin{aligned}
\omega&=e^{12}+e^{34}+e^{56},\\
\psi_+&=e^{135}-e^{146}-e^{236}-e^{245},\\
\psi_-&=e^{136}+e^{145}+e^{235}-e^{246}.
\end{aligned}
\end{equation*}

\begin{proposition}\label{G2cerrada}
Let $(\mathfrak{h}, \omega, \psi_{+})$ be a symplectic half-flat Lie algebra, and let $D$ be 
a derivation of $\mathfrak{h}$ such that $D$ is the real representation of $A \in \mathfrak{sl}(3, \mathbb{C}), $ with respect to a basis  $\{e_1,\dots, e_6\}$ of $\frak h$ such that     $\omega, \, \psi_ +$ and $\psi_-$ have the canonical expression. Then, the Lie algebra
\begin{equation*}
\mathfrak{g}=\mathfrak{h}\oplus_D\mathbb{R}\xi,
\end{equation*}
has a closed $\mathrm{G}_2$ form.

\begin{proof}
By \eqref{sl(3,C)} we have
$$
 D = \left(
\begin{array}{cc|cc|cc}
 a_{1,1} & a_{1,2} & a_{1,3} & a_{1,4} & a_{1,5} & a_{1,6}  \\
 -a_{1,2} & a_{1,1} & -a_{1,4} & a_{1,3} & -a_{1,6} & a_{1,5}  \\\hline
 a_{3,1} & a_{3,2} & a_{3,3} & a_{3,4} & a_{3,5} & a_{3,6}  \\
 -a_{3,2} & a_{3,1} & -a_{3,4} & a_{3,3} & -a_{3,6} & a_{3,5}  \\\hline
 a_{5,1} & a_{5,2} & a_{5,3} & a_{5,4} & -a_{1,1}-a_{3,3} & -a_{1,2}-a_{3,4}  \\
 -a_{5,2} & a_{5,1} & -a_{5,4} & a_{5,3} & a_{1,2}+a_{3,4} & -a_{1,1}-a_{3,3}  \\
\end{array}
\right),
$$
with respect to the basis  $\{e_1,\dots, e_6\}$ of $\frak h$ such that the $\mathrm{SU}(3)$-structure ($\omega, \psi_ +$) has the canonical expression.

Consider on $\mathfrak{g}$, the $\mathrm{G}_2$  form
\begin{equation}\label{G2-form-algebra}
\varphi=\omega\wedge \eta +\psi_+,
\end{equation}
with $\eta$ the 1-form such that $\eta(X)=0$ for all $X\in \mathfrak{h}$ and $\eta(\xi)=1$.\\
For every $U,V,W,T\in \mathfrak{h}$
\begin{equation*}
d\varphi(U,V,W,T)=d\psi_+(U,V,W,T),
\end{equation*}
which vanishes since $\psi_+$ is closed.\\
Hence, consider
\begin{equation*}
\begin{aligned}
d\varphi(U,V,W,\xi)=&-\varphi([U,V],W,\xi)+\varphi([U,W],V,\xi)-\varphi([U,\xi],V,W)\\
&-\varphi([V,W],U,\xi)+\varphi([V,\xi],U,W)-\varphi([W,\xi],U,V),
\end{aligned}
\end{equation*}
which by definition of $\varphi$ is exactly
\begin{equation*}
\begin{aligned}
&-\omega([U,V],W)+\omega([U,W],V)-\omega([V,W],U)-\psi_+([U,\xi],V,W)\\
&+\psi_+([V,\xi],U,W)-\psi_+([W,\xi],U,V)=d\omega(U,V,W)+\psi_+(D(U),V,W)\\
&+\psi_+(U,D(V),W)+\psi_+(U,V,D(W)).
\end{aligned}
\end{equation*}
Therefore, since $\omega$ is closed
\begin{equation*}
d\varphi(U,V,W,\xi)=\psi_+(D(U),V,W)+\psi_+(U,D(V),W)+\psi_+(U,V,D(W)).
\end{equation*}
Taking into account the expression of  $D$, and $\psi_+$ in terms of the basis $\{e_1,\dots, e_6\}$ an easy computation shows that
\begin{equation*}
\psi_+(D(e_i),e_j,e_k)+\psi_+(e_i,D(e_j),e_k)+\psi_+(e_i,e_j,D(e_k))=0,
\end{equation*}
for every triple $(e_i, e_j, e_k)$ of elements of the basis of $\mathfrak{h}$. Thus, the $\mathrm{G}_2$  form $\varphi$ defined in \eqref{G2-form-algebra} is closed in $\mathfrak{g}$.
\end{proof}
\end{proposition}

The previous proposition describes a method to construct new 7-dimensional Lie algebras with a closed $\mathrm{G}_2$-structure.

\begin{remark}
Note that the trace of $D$, the real representation of certain $A \in \mathfrak{sl}(3, \mathbb{C})$ vanishes. Therefore, the Lie algebra $\frg=\frh \oplus_{D} \mathbb{R}e_7$ will be unimodular if and only if $\frh$ is so.
\end{remark}

\begin{example}\label{ejemplo1}
Next, we show a {\em new example of compact
solvmanifold
with closed $\mathrm{G}_2$ form}.
Let $\mathfrak{h}$ 
be the 6-dimensional abelian Lie algebra defined by the structure equations
\begin{equation*}
\mathfrak{h}=(0, 0, 0, 0, 0, 0),
\end{equation*}
The almost Hermitian structure $(g, J)$ on $\mathfrak{h}$ given by
\begin{equation*}
g=\sum_{i=1}^6 e^i\otimes e^i,  \qquad  Je_1=e_2, \qquad Je_3=e_4, \qquad Je_5=e_6
\end{equation*}
is such that its K\"ahler form is 
\begin{equation*}
\omega=e^{12}+e^{34}+e^{56}.
\end{equation*}
Thus, $(g,J)$ together with the complex volume form $\Psi=\psi_{+} +i\,\psi_{-}$, where
\begin{equation*}
\begin{aligned}
\psi_+&=e^{135}-e^{146}-e^{236}-e^{245},\\
\psi_-&=e^{136}+e^{145}+e^{235}-e^{246},
\end{aligned}
\end{equation*}
define an $\mathrm{SU}(3)$-structure on $\mathfrak{h}$. Clearly, $d\,\omega^2=d\,\psi_+=0$, so $(g, J, \Psi=\psi_{+}+i\, \psi_-)$
is a half-flat $\mathrm{SU}(3)$-structure on $\mathfrak{h}$.
\\
Consider now the derivation $D$ of $\mathfrak{h}$ given by
 \begin{equation*}
\left(
\begin{array}{cc|cc|cc}
 1&&&&&  \\
 &1&&&& \\\hline
 &&-1&&&  \\
&&&-1&&  \\\hline
 &&&&0&  \\
 &&&&&0  \\
\end{array}
\right) \in \mathfrak{sl}(3, \mathbb{C}),
\end{equation*}
that is,
\begin{equation*}
D(e_1)=e_1,  \quad D(e_2)=e_2,  \quad D(e_3)=-e_3, \quad D(e_4)=-e_4,
\end{equation*}
Take the Lie algebra 
\begin{equation*}
\mathfrak{g}=\mathfrak{h}\oplus_D \mathbb{R}e_7,
\end{equation*}
whose structure equations are
\begin{equation*}
\mathfrak{g}=(e^{17}, e^{27},-e^{37}, -e^{47},0,0,0).
\end{equation*}
Then, the 3-form $\varphi$ given by
\begin{equation*}
\varphi=e^{127}+e^{347}+e^{567}+e^{136}+e^{145}+e^{235}-e^{246}
\end{equation*}
is a closed $\mathrm{G}_2$  form on $\mathfrak{g}$. 

Lets denote by $G$ the simply connected and completely solvable Lie group consisting on matrices of the form.
\begin{equation*}
a=\left(
\begin{array}{cc|cc|cc|cc}
e^{x_7} &  &  &  &  & &&x_1  \\
  & e^{x_7} & &  &  & &&x_2  \\\hline
 &  & e^{-x_7} &  &  & &&x_3  \\
  &  &  & e^{-x_7} &  & &&x_4 \\\hline
  &  &  & & 1 &  &&x_5 \\
  &  &  &  & & 1  &&x_6\\\hline
    &  &  &  & &  &1&x_7\\
      &  &  &  & &   &&1\\
\end{array}
\right), 
\end{equation*}
with $x_i \in \mathbb{R}$, for $i=1,\dots, 7$. Then a global system of coordinates $\{x_i\}$ for $G$ is defined by $x_i(a)=x_i$. An standard calculation shows that a basis for the left invariant 1-forms on $G$ can be described by
\begin{align*}
&e^1=e^{-x_7}dx_1, \quad e^2=e^{-x_7}dx_2, \quad e^3=e^{x_7}dx_3, \quad e^4=e^{x_7}dx_4,\\
&e^5=dx_5, \qquad e^6=dx_6, \qquad  and \qquad e^7=dx_7.
\end{align*}
Therefore $\mathfrak{g}$ is exactly the Lie algebra of $G$.
Notice that $G= \mathbb{R} \ltimes_{\phi} \mathbb{R}^6$, where $\mathbb{R}$ acts on $\mathbb{R}^6$ via $\phi_t$ described by 
 \begin{equation*}
\phi_t=
\left(
\begin{array}{cc|cc|cc}
e^{t} &  &  &  &  &   \\
  & e^{t} & &  &  &   \\\hline
 &  & e^{-t} &  &  &   \\
  &  &  & e^{-t} &  &  \\\hline
  &  &  & & 1 &   \\
  &  &  &  & & 1  \\
\end{array}
\right). 
\end{equation*}
Thus the operation on the group $G$ is given by 
\begin{equation*}
a \cdot b= (b_1e^{a_7}+a_1, b_2e^{a_7}+a_2, b_3e^{-a_7}+a_3, b_4e^{-a_7}+a_4, b_5+a_5, b_6+a_6, b_7+a_7),
\end{equation*}
where $a=(a_1,\dots, a_7)$ and $b=(b_1,\dots, b_7)$.\\ 
To construct the lattice $\Gamma$ of $G$ it is enough to find some real number $t_0$ such that $\phi_{t_0}$ is conjugated to an element $A\in SL(6,\mathbb{Z})$. If $\Gamma_0$ denotes a lattice of $\mathbb{R}^6$ invariant under $\phi_{t_0}$, take
\begin{equation*}
\Gamma=(t_0 \, \mathbb{Z}) \ltimes_{\phi} \Gamma_0.
\end{equation*}
Consider the matrix
\begin{equation*}
A=
\left(
\begin{array}{cc|cc|cc}
2 & 1 &  &  &  &   \\
 1 & 1 & &  &  &   \\\hline
 &  & 2 & 1 &  &   \\
  &  & 1 & 1&  &  \\\hline
  &  &  & & 1 &   \\
  &  &  &  &  & 1  \\
\end{array}
\right), 
\end{equation*}
with double eigenvalues $\frac{3+\sqrt{5}}{2}, \frac{3-\sqrt{5}}{2}$. Taking $t_0=Ln(\frac{3+\sqrt{5}}{2})$ we have that $e^{t_0D}$ and $A$ are conjugated. In particular, take
\begin{equation*}
P=
\left(
\begin{array}{cc|cc|cc}
1 & \frac{-1+\sqrt{5}}{2} &  &  &  &   \\
 1 & \frac{-1-\sqrt{5}}{2} & &  &  &   \\\hline
 &  & 1 & \frac{-1+\sqrt{5}}{2} &  &   \\
  &  & 1 & \frac{-1-\sqrt{5}}{2} &  &  \\\hline
  &  &  & & 1 &  \\
  &  &  &  &  & 1  \\
\end{array}
\right), 
\end{equation*}
it is easy to check that $PA=\phi_{t_0}P$. So, the lattice defined by
\begin{equation*}
\Gamma_0=P \, \mathbb{Z}\langle e_1,\dots, e_6 \rangle
\end{equation*}
is invariant under the group $t_0\mathbb{Z}$. Thus 
\begin{equation*}
\Gamma=(t_0 \, \mathbb{Z})\ltimes_{\phi}\Gamma_0
\end{equation*}
is a lattice of $G$. Then, the compact solvmanifold $S=\Gamma \backslash G$ admits a closed $\mathrm{G}_2$-structure.
\end{example}

\medskip

Regarding the converse of Proposition \ref{G2cerrada} we have the following result

\begin{proposition}\label{G2cerradainversa}
Let $\varphi$ be the closed $\mathrm{G}_2$ form
\begin{equation*}
\varphi=e^{127}+e^{347}+e^{567}+e^{135}-e^{146}-e^{236}-e^{245},
\end{equation*}
on the 7-dimensional Lie algebra
\begin{equation*}
\mathfrak{g}=\mathfrak{h}\oplus_D \mathbb{R}e_7,
\end{equation*}
where $D$ is a derivation of $\mathfrak{h}=\frg/\langle e_7 \rangle$ such that it is given with respect to the basis $\{e_1, \dots, e_6\}$ by a matrix of the form \eqref{sl(3,C)}. Then, the 6-dimensional Lie algebra $\frh$ has a symplectic half-flat structure.
\begin{proof}
Consider the pair of forms $(\omega, \psi_+)$ on $\frh$ defined by 
\begin{equation*}
\omega=\iota_{e_7}\varphi \qquad \text{ and } \qquad \psi_+=\pi^* \varphi,
\end{equation*}
with $\pi$ the projection of $\frg$ onto $\frh$. 
Thus $(\omega, \psi_+)$ define an $\mathrm{SU}(3)$-structure on $\frh$, which in terms of the basis $\{e_1, \dots, e_6\}$ of $\frh$ has the canonical expression, that is 
\begin{equation*}
\begin{aligned}
\omega&=e^{12}+e^{34}+e^{56},\\
\psi_+&=e^{135}-e^{146}-e^{236}-e^{245}.
\end{aligned}
\end{equation*}
Then, for every $U,V,W,T\in \mathfrak{h}$
\begin{equation*}
d\psi_+(U,V,W,T)=d\varphi(U,V,W,T),
\end{equation*}
which vanishes since $\varphi$ is closed, and therefore $\psi_+$ is closed.

Taking into account the expression of $D$, and $\psi_+$ in terms of the basis $\{e_1,\dots, e_6\}$ an easy computation shows that
\begin{equation*}
\psi_+(D(e_i),e_j,e_k)+\psi_+(e_i,D(e_j),e_k)+\psi_+(e_i,e_j,D(e_k))=0,
\end{equation*}
for every triple $(e_i, e_j, e_k)$ of elements of the basis of $\mathfrak{h}$. Consider now $d\varphi(U,V,W, e_7)$ which is exactly
\begin{equation*}
\begin{aligned}
d\varphi(U,V,W,e_7)=&d \, \omega(U,V,W)+\psi_+(D(U),V,W)\\
&+\psi_+(U,D(V),W)+\psi_+(U,V,D(W)).
\end{aligned}
\end{equation*}
Therefore, since $D$ is the real representation of certain $A \in \mathfrak{sl}(3,\mathbb{C})$ 
\begin{equation*}
d\omega(U,V,W)=d\varphi(U,V,W,e_7),
\end{equation*}
which vanishes since $\varphi$ is closed. Thus the $\mathrm{SU}(3)$-structure on $\mathfrak{h}$ given by the pair $(\omega, \psi_+)$ is symplectic half-flat.
\end{proof}
\end{proposition}

\begin{example}\label{ejemplo2}
Let $\frg$ be the 7-dimensional nilpotent Lie algebra described by the structure equations 
\begin{equation*}
\frg=(0,0,e^{17},e^{15}+e^{27},0,e^{13},0). 
\end{equation*}
Then, the $\mathrm{G}_2$ form given by 
\begin{equation*}
\varphi=e^{127}+e^{347}+e^{567}+e^{135}-e^{146}-e^{236}-e^{245},
\end{equation*}
is closed.\\
Notice that $\frg$ is of the form
\begin{equation*}
\frg=\frh\oplus_D\mathbb{R}e_7,
\end{equation*}
where the 6-dimensional Lie algebra $\frh=\frg/\langle e_7 \rangle$ is described by the structure equations 
\begin{equation*}
\frh=(0,0,0,e^{15},0,e^{13}). 
\end{equation*}
The derivation $D$ of $\frh$, is given with respect to the basis $\{e_1,\dots,e_6\}$ by the matrix
\begin{equation*}
D=\left(
\begin{array}{cc|cc|cc}
 &&1&&&  \\
 &&&1&& \\\hline
 &&&&&  \\
&&&&&  \\\hline
 &&&&&  \\
 &&&&&  \\
\end{array}
\right).
\end{equation*}
Therefore $D$ is the real representation of certain $A \in \mathfrak{sl}(3,\mathbb{C})$. Then, the $\mathrm{SU}(3)$-structure defined by
\begin{equation*}
\begin{aligned}
\omega&=e^{12}+e^{34}+e^{56},\\
\psi_+&=e^{135}-e^{146}-e^{236}-e^{245}.
\end{aligned}
\end{equation*}
is symplectic half-flat.
\end{example}

As a consequence of Proposition \ref{G2cerrada} and \ref{G2cerradainversa} we have the following result:

\begin{theorem}
Let $\frh$ be a 6-dimensional Lie algebra and let $\frg$ be a 7-dimensional Lie algebra satisfying 
\begin{equation*}
\frg=\frh\oplus_D\mathbb{R}e_7,
\end{equation*}
with $D$ a derivation of $\mathfrak{h}$ given by \eqref{sl(3,C)} in terms of a basis $\{e_1,\dots,e_6\}$ of $\frh$. Then the following two conditions are equivalent:
\begin{enumerate}
\item The $\mathrm{SU}(3)$-structure on $\frh$ given by 
\begin{equation*}
\begin{aligned}
\omega&=e^{12}+e^{34}+e^{56},\\
\psi_+&=e^{135}-e^{146}-e^{236}-e^{245}.
\end{aligned}
\end{equation*}
is symplectic half-flat.
\item The $\mathrm{G}_2$-structure on $\frg$ given by 
\begin{equation*}
\varphi=e^{127}+e^{347}+e^{567}+e^{135}-e^{146}-e^{236}-e^{245},
\end{equation*}
is closed.
\end{enumerate}
\end{theorem}

\end{section}


\begin{section}{Lie algebras with a coclosed $\mathrm{G}_2$-structure} \label{Coclosed}
In this section we show that if a 6-dimensional half-flat Lie algebra is endowed with a particular type of derivation, then a 
Lie algebra  with a coclosed $\mathrm{G}_2$-structure can be constructed.

We recall that a coclosed $\mathrm{G}_2$-structure on a real Lie algebra $\frg$ of dimension 7 consists on the presence of a $\mathrm{G}_2$ form
which is coclosed. In order to obtain an expression adapted to our purposes, in this section we charaterize a $\mathrm{G}_2$ form on $\frg$ as 
a 3-form that can be written as
\begin{equation*}
\varphi=e^{127}+e^{347}+e^{567}+e^{246}-e^{235}-e^{145}-e^{136},
\end{equation*}
with respect to some basis $\{e^1,\dots,e^7\}$ of the dual space of $\frg$.

\begin{proposition}\label{G2cocerrada}
Let $(\mathfrak{h}, \omega, \psi_{+})$ be a half-flat Lie algebra, and let $D$ be 
a derivation of $\mathfrak{h}$ such that $D \in \mathfrak{sp}(6, \mathbb{R})$,  with respect to a basis  $\{e_1,\dots, e_6\}$ of $\frak h$ satisfying
 that $\omega, \, \psi_ + $ and $\psi_-$ have the canonical expression. Then, the Lie algebra
\begin{equation*}
\mathfrak{g}=\mathfrak{h}\oplus_D\mathbb{R}\xi,
\end{equation*}
has a coclosed $\mathrm{G}_2$ form.

\begin{proof}
Since $D\in \mathfrak{sp}(6,\mathbb{R})$ with respect  to the  basis $\{e_1, \ldots, e_6\}$, we can write 
\begin{equation}\label{sp(6,R)}
D=
\left(
\begin{array}{cc|cc|cc}
 a_{1,1} & a_{1,2} & a_{1,3} & a_{1,4} & a_{1,5} & a_{1,6}  \\
 a_{2,1} & -a_{1,1} & a_{2,3} & a_{2,4} & a_{2,5} & a_{2,6}  \\\hline
 -a_{2,4} & a_{1,4} & a_{3,3} & a_{3,4} & a_{3,5} & a_{3,6}  \\
 a_{2,3} & -a_{1,3} & a_{4,3} & -a_{3,3} & a_{4,5} & a_{4,6}  \\\hline
 -a_{2,6} & a_{1,6} & -a_{4,6} & a_{3,6} & a_{5,5} & a_{5,6}  \\
 a_{2,5} & -a_{1,5} & a_{4,5} & -a_{3,5} & a_{6,5} & -a_{5,5}  \\
\end{array}
\right).
\end{equation}
Consider on $\mathfrak{g}$, the $\mathrm{G}_2$  form 
\begin{equation}\label{G2-form-algebra2}
\varphi=\omega\wedge \eta- \psi_-,
\end{equation}
thus
\begin{equation*}
\ast\varphi=\frac{1}{2}\omega\wedge\omega +\psi_+\wedge \eta, 
\end{equation*}
where $\eta$ is the 1-form satisfying that $\eta(X)=0$ for all $X\in \mathfrak{h}$ and $\eta(\xi)=1$.\\
For every $U,V,W,T,R \in \mathfrak{h}$, 
\begin{equation*}
d\ast\varphi(U,V,W,T,R)= d \, \omega\wedge\omega(U,V,W,T,R),
\end{equation*}
which vanishes since $\omega\wedge \omega$ is closed.\\
Hence, consider
\begin{equation*}
\begin{aligned}
d\ast\varphi(U,V,W,T,\xi)=&-\ast\varphi([U,V],W,T,\xi)+\ast\varphi([U,W],V,T,\xi)-\ast\varphi([U,T],V,W,\xi)\\
&+\ast\varphi([U,\xi],V,W,T)-\ast\varphi([V,W],U,T,\xi)+\ast\varphi([V,T],U,W,\xi)\\
&-\ast\varphi([V,\xi],U,W,T)-\ast\varphi([W,T],U,V,\xi)+\ast\varphi([W,\xi],U,V,T)\\
&-\ast\varphi([T,\xi],U,V,W),
\end{aligned}
\end{equation*}
which by the definition of $\ast\varphi$ is exactly
\begin{equation*}
\begin{aligned}
&-\psi_+([U,V],W,T)+\psi_+([U,W],V,T)-\psi_+([U,T],V,W)-\psi_+([V,W],U,T)\\
&+\psi_+([V,T],U,W)-\psi_+([W,T],U,V)+\frac{1}{2}\omega\wedge\omega([U,\xi],V,W,T)\\
&-\frac{1}{2}\omega\wedge\omega([V,\xi],U,W,T)+\frac{1}{2}\omega\wedge\omega([W,\xi],U,V,T)-\frac{1}{2}\omega\wedge\omega([T,\xi],U,V,W)\\
&=d\psi_+(U,V,W,T)+\frac{1}{2}\omega\wedge\omega(D(U),V,W,T)+\frac{1}{2}\omega\wedge\omega(U,D(V),W,T)\\
&+\frac{1}{2}\omega\wedge\omega(U,V,D(W),T)+\frac{1}{2}\omega\wedge\omega(U,V,W,D(T)).
\end{aligned}
\end{equation*}
Therefore since $\psi_+$ is closed 
\begin{equation*}
\begin{aligned}
d\ast\varphi(U,V,W,T,\xi)=&\frac{1}{2}\omega\wedge\omega(D(U),V,W,T)+\frac{1}{2}\omega\wedge\omega(U,D(V),W,T)\\
&+\frac{1}{2}\omega\wedge\omega(U,V,D(W),T)+\frac{1}{2}\omega\wedge\omega(U,V,W,D(T)).
\end{aligned}
\end{equation*}
Using the expressions of $D$ and $\omega$ with respect to the basis $\{ e_1, \ldots, e_6 \}$, can be checked that
\begin{equation*}
\begin{aligned}
&\omega\wedge\omega(D(e_i),e_j,e_k,e_l)+\omega\wedge\omega(e_i,D(e_j),e_k,e_l)\\
&+\omega\wedge\omega(e_i,e_j,D(e_k),e_l)+\omega\wedge\omega(e_i,e_j,e_k,D(e_l))=0,
\end{aligned}
\end{equation*}
for every quadruplet $(e_i, e_j, e_k, e_l)$ of elements of the basis of $\mathfrak{h}$. Thus, the $\mathrm{G}_2$  form $\varphi$ 
defined in \eqref{G2-form-algebra2} is coclosed in $\mathfrak{g}$.
\end{proof}
\end{proposition}

Previous proposition describes a method to construct 7-dimensional Lie algebras with a coclosed $\mathrm{G}_2$-structure.

\begin{remark}
Note that the trace of $D \in \mathfrak{sp}(6, \mathbb{R})$ vanishes. Therefore, the Lie algebra $\frg=\frh \oplus_D \mathbb{R}e_7$ will be unimodular if and only if $\frh$ is so. 
\end{remark}
  
\begin{example}\label{ejemplo3}
Next, we show a {\em new example of 
solvable Lie algebra
with coclosed $\mathrm{G}_2$ form}.
Let $\mathfrak{h}$ 
be the 6-dimensional abelian Lie algebra defined by the structure equations
\begin{equation*}
\mathfrak{h}=(0, 0, 0, 0, 0, 0),
\end{equation*}
The almost Hermitian structure $(g, J)$ on $\mathfrak{h}$ given by
\begin{equation*}
g=\sum_{i=1}^6 e^i\otimes e^i,  \qquad  Je_1=e_2, \qquad Je_3=e_4, \qquad Je_5=e_6
\end{equation*}
is such that its K\"ahler form is 
\begin{equation*}
\omega=e^{12}+e^{34}+e^{56}.
\end{equation*}
Thus, $(g,J)$ together with the complex volume form $\Psi=\psi_{+} +i\,\psi_{-}$, where
\begin{equation*}
\begin{aligned}
\psi_+&=e^{135}-e^{146}-e^{236}-e^{245},\\
\psi_-&=e^{136}+e^{145}+e^{235}-e^{246},
\end{aligned}
\end{equation*}
define an $\mathrm{SU}(3)$-structure on $\mathfrak{h}$. Clearly, $d\,\omega^2=d\,\psi_+=0$, so $(g, J, \Psi=\psi_{+}+i\, \psi_-)$
is a half-flat $\mathrm{SU}(3)$-structure on $\mathfrak{h}$.
\\
Consider now the derivation $D$ of $\mathfrak{h}$ given by
 \begin{equation*}
\left(
\begin{array}{cc|cc|cc}
 1&&&&&  \\
 &-1&&&& \\\hline
 &&1&&&  \\
&&&-1&&  \\\hline
 &&&&1&  \\
 &&&&&-1  \\
\end{array}
\right) \in \mathfrak{sp}(6, \mathbb{R}),
\end{equation*}
that is,
\begin{equation*}
\begin{aligned}
&D(e_1)=e_1,  \qquad D(e_2)=-e_2,  \qquad D(e_3)=e_3, \\ 
&D(e_4)=-e_4,   \quad D(e_5)=e_5,  \quad D(e_6)=-e_6.
\end{aligned}
\end{equation*}
Take the Lie algebra 
\begin{equation*}
\mathfrak{g}=\mathfrak{h}\oplus_D \mathbb{R}e_7,
\end{equation*}
whose structure equations are
\begin{equation*}
\mathfrak{g}=(e^{17}, -e^{27},e^{37}, -e^{47},e^{57},-e^{67},0).
\end{equation*}
Then, the 3-form $\varphi$ given by
\begin{equation*}
\varphi=e^{127}+e^{347}+e^{567}+e^{136}+e^{145}+e^{235}-e^{246}
\end{equation*}
is a closed $\mathrm{G}_2$  form on $\mathfrak{g}$. 
\end{example}

\begin{example}\label{ejemplo3}
Let $\mathfrak{h}$  be the 6-dimensional abelian Lie algebra described by the structure equations
\begin{equation*}
\mathfrak{h}=(0,0, 0, 0,0,0).
\end{equation*}
The almost Hermitian structure given by
\begin{equation*}
g=\sum_{i=1}^6 e^i\otimes e^i \quad \text{and} \quad Je_1=e_2, \quad Je_3=e_4, \quad Je_5=e_6,
\end{equation*}
is such that its K\"ahler form is 
\begin{equation*}
\omega=e^{12}+e^{34}+e^{56}.
\end{equation*}
Thus, $(g,J)$ together with the complex volume form $\Psi=\psi_+ +i \, \psi_-$ where
\begin{equation*}
\begin{aligned}
\psi_+&=e^{135}-e^{146}-e^{236}-e^{245},\\
\psi_-&=e^{136}+e^{145}+e^{235}-e^{246},
\end{aligned}
\end{equation*}
define an $\mathrm{SU}(3)$-structure on $\mathfrak{h}$. Concretely, since $\omega^2$ and $\psi_+$ are closed it is a half-flat structure. 
Consider now the derivation $D$ of $\mathfrak{h}$ given by
\begin{equation*}
D=\left(
\begin{array}{cc|cc|cc}
 1&&&&&  \\
 &-1&&&& \\\hline
 &&1&&&  \\
&&&-1&&  \\\hline
 &&&&1&  \\
 &&&&&-1  \\
\end{array}
\right) \in \mathfrak{sp}(6, \mathbb{R}),
\end{equation*}
that is,
\begin{equation*}
\begin{aligned}
&D(e_1)=e_1, \qquad  D(e_2)=-e_2, \qquad D(e_3)=e_3, \\
&D(e_4)=-e_4, \quad  D(e_5)=e_5 \quad  and \quad  D(e_6)=-e_6.
\end{aligned}
\end{equation*}
Thus, the Lie algebra 
\begin{equation*}
\mathfrak{g}=\mathfrak{h}\oplus_D \mathbb{R}e_7,
\end{equation*}
which is described by the structure equations
\begin{equation*}
\mathfrak{g}=(e^{17}, -e^{27},e^{37}, -e^{47},e^{57},-e^{67},0),
\end{equation*}
is completely solvable and admits the coclosed $\mathrm{G}_2$  form
\begin{equation*}
\varphi=e^{127}+e^{347}+e^{567}+e^{136}+e^{145}+e^{235}-e^{246}.
\end{equation*}
Lets denote by $G$ the simply connected and completely solvable Lie group consisting on matrices of the form.
\begin{equation*}
a=\left(
\begin{array}{cc|cc|cc|cc}
e^{x_7} &  &  &  &  & &&x_1  \\
  & e^{-x_7} & &  &  & &&x_2  \\\hline
 &  & e^{x_7} &  &  & &&x_3  \\
  &  &  & e^{-x_7} &  & &&x_4 \\\hline
  &  &  & & e^{x_7} &  &&x_5 \\
  &  &  &  & & e^{-x_7}  &&x_6\\\hline
    &  &  &  & &  &1&x_7\\
      &  &  &  & &   &&1\\
\end{array}
\right), 
\end{equation*}
with $x_i \in \mathbb{R}$, for $i=1,\dots, 7$. Then a global system of coordinates $\{x_i\}$ for $G$ is defined by $x_i(a)=x_i$. An standard calculation shows that a basis for the left invariant 1-forms on $G$ can be described by
\begin{align*}
&e^1=e^{-x_7}dx_1, \quad e^2=e^{x_7}dx_2, \quad e^3=e^{-x_7}dx_3, \quad e^4=e^{x_7}dx_4,\\
&e^5=e^{-x_7}dx_5, \quad e^6=e^{-x_7}dx_6, \qquad  and \qquad e^7=dx_7.
\end{align*}
Therefore $\mathfrak{g}$ is exactly the Lie algebra of $G$.
Notice that $G= \mathbb{R} \ltimes_{\phi} \mathbb{R}^6$, where $\mathbb{R}$ acts on $\mathbb{R}^6$ via $\phi_t$ described by 
 \begin{equation*}
\phi_t=
\left(
\begin{array}{cc|cc|cc}
e^{t} &  &  &  &  &   \\
  & e^{-t} & &  &  &   \\\hline
 &  & e^{t} &  &  &   \\
  &  &  & e^{-t} &  &  \\\hline
  &  &  & & e^{t} &   \\
  &  &  &  & & e^{-t}  \\
\end{array}
\right). 
\end{equation*}
Thus the operation on the group $G$ is given by 
\begin{equation*}
a \cdot b= (b_1e^{a_7}+a_1, b_2e^{-a_7}+a_2, b_3e^{a_7}+a_3, b_4e^{-a_7}+a_4, b_5e^{a_7}+a_5, b_6e^{-a_7}+a_6, b_7+a_7),
\end{equation*}
where $a=(a_1,\dots, a_7)$ and $b=(b_1,\dots, b_7)$.\\ 
To construct the lattice $\Gamma$ of $G$ it is enough to find some real number $t_0$ such that $\phi_{t_0}$ is conjugated to an element $A\in SL(6,\mathbb{Z})$. If $\Gamma_0$ denotes a lattice of $\mathbb{R}^6$ invariant under $\phi_{t_0}$, take
\begin{equation*}
\Gamma=(t_0 \, \mathbb{Z}) \ltimes_{\phi} \Gamma_0.
\end{equation*}
Consider the matrix
\begin{equation*}
A=
\left(
\begin{array}{cc|cc|cc}
2 & 1 &  &  &  &   \\
 1 & 1 & &  &  &   \\\hline
 &  & 2 & 1 &  &   \\
  &  & 1 & 1&  &  \\\hline
  &  &  & & 2 &  1 \\
  &  &  &  & 1 & 1  \\
\end{array}
\right), 
\end{equation*}
with triple eigenvalues $\frac{3+\sqrt{5}}{2}, \frac{3-\sqrt{5}}{2}$. Taking $t_0=Ln(\frac{3+\sqrt{5}}{2})$ we have that $e^{t_0D}$ and $A$ are conjugated. In particular, take
\begin{equation*}
P=
\left(
\begin{array}{cc|cc|cc}
1 & \frac{-1+\sqrt{5}}{2} &  &  &  &   \\
 1 & \frac{-1-\sqrt{5}}{2} & &  &  &   \\\hline
 &  & 1 & \frac{-1+\sqrt{5}}{2} &  &   \\
  &  & 1 & \frac{-1-\sqrt{5}}{2} &  &  \\\hline
  &  &  & & 1 &  \frac{-1+\sqrt{5}}{2} \\
  &  &  &  & 1 & \frac{-1-\sqrt{5}}{2}  \\
\end{array}
\right), 
\end{equation*}
it is easy to check that $PA=\phi_{t_0}P$. So, the lattice defined by
\begin{equation*}
\Gamma_0=P \, \mathbb{Z}\langle e_1,\dots, e_6 \rangle
\end{equation*}
is invariant under the group $t_0\mathbb{Z}$. Thus 
\begin{equation*}
\Gamma=(t_0 \, \mathbb{Z})\ltimes_{\phi}\Gamma_0
\end{equation*}
is a lattice of $G$. Then, the compact solvmanifold $S=\Gamma \backslash G$ admits a coclosed $\mathrm{G}_2$-structure.
\end{example}

Considering the converse of Proposition \ref{G2cocerrada} we obtain the next result.
\begin{proposition}\label{G2cocerradainversa}
Let $\varphi$ be the coclosed $\mathrm{G}_2$ form
\begin{equation*}
\varphi=e^{127}+e^{347}+e^{567}+e^{136}-e^{145}-e^{235}-e^{246},
\end{equation*}
on the 7-dimensional Lie algebra
\begin{equation*}
\mathfrak{g}=\mathfrak{h}\oplus_D \mathbb{R}e_7,
\end{equation*}
where $D$ is a derivation of $\mathfrak{h}=\frg/\langle e_7 \rangle$ such that it is given with respect to the basis $\{e_1, \dots, e_6\}$ by a matrix of the form \eqref{sp(6,R)}. 
Then, the 6-dimensional Lie algebra $\frh$ has a half-flat structure.
\begin{proof}
Consider the pair of forms $(\omega, \psi_+)$ on $\frh$ defined by 
\begin{equation*}
\omega=\iota_{e_7}\varphi \qquad \text{ and } \qquad \psi_+=\ast \psi_-,
\end{equation*}
where
\begin{equation*}
\psi_-=-\pi^* \varphi,
\end{equation*}
with $\pi$ the projection of $\frg$ onto $\frh$. 
Thus $(\omega, \psi_+)$ define an $\mathrm{SU}(3)$-structure on $\frh$, which in terms of the basis $\{e_1, \dots, e_6\}$ of $\frh$ has the canonical expression, that is 
\begin{equation*}
\begin{aligned}
\omega&=e^{12}+e^{34}+e^{56},\\
\psi_+&=e^{135}-e^{146}-e^{236}-e^{245}.
\end{aligned}
\end{equation*}
Therefore we have that
\begin{equation*}
\ast\varphi=\frac12\omega\wedge\omega+\psi_+\wedge e^7,
\end{equation*}
and then, for every $U,V,W,T,R \in \mathfrak{h}$
\begin{equation*}
d\omega\wedge\omega(U,V,W,T,R)=d\ast\varphi(U,V,W,T,R),
\end{equation*}
which vanishes since $\ast\varphi$ is closed. Thus, $\omega \wedge \omega$ is closed.

Consider now $d\ast\varphi(U,V,W,T,e_7)$ which by definition of $\ast \varphi$ is exactly
\begin{equation*}
\begin{aligned}
d\ast\varphi(U,V,W,T,e_7)=& \, \omega\wedge \omega(D(U),V,W,T)+\omega\wedge \omega(U,D(V),W,T)\\
&+\omega\wedge \omega(U,V,D(W),T)+\omega\wedge \omega(U,V,W,D(T))\\
&+ d\psi_+(U,V,W,T).
\end{aligned}
\end{equation*}
Since $D\in \mathfrak{sp}(6,\mathbb{R})$
\begin{equation*}
d\psi_+(U,V,W,T)=d\ast\varphi(U,V,W,T,e_7),
\end{equation*}
that vanishes because $\varphi$ is coclosed. Therefore the $\mathrm{SU}(3)$-structure on $\frh$ given by the pair $(\omega, \psi_+)$ is half-flat.
\end{proof}
\end{proposition}

\begin{example}\label{ejemplo4}
Let $\frg$ be the 7-dimensional Lie algebra described by the structure equations 
\begin{equation*}
\frg=(e^{35}+e^{46},0,e^{67},e^{57},e^{47},e^{37},0), 
\end{equation*}
then, the $\mathrm{G}_2$ form given by 
\begin{equation*}
\varphi=e^{127}+e^{347}+e^{567}+e^{246}-e^{235}-e^{136}-e^{145},
\end{equation*}
is coclosed, that is, the 4-form
\begin{equation*}
\ast \varphi=e^{1234}+e^{1256}+e^{3456}+e^{1357}-e^{1467}-e^{2367}-e^{2457},
\end{equation*}
is closed.
Notice that $\frg$ is of the form
\begin{equation*}
\frg=\frh\oplus_D\mathbb{R}e_7,
\end{equation*}
where the 6-dimensional Lie algebra $\frh=\frg/\langle e_7 \rangle$ is described by the structure equations 
\begin{equation*}
\frh=(e^{35}+e^{46},0,0,0,0,0). 
\end{equation*}
The derivation $D$ of $\frh$, is given with respect to the basis $\{e_1,\dots,e_6\}$ by the matrix
\begin{equation*}
D=\left(
\begin{array}{cc|cc|cc}
 &&&&&  \\
 &&&&& \\\hline
 &&&&&1  \\
&&&&1&  \\\hline
 &&&1&&  \\
 &&1&&&  \\
\end{array}
\right).
\end{equation*}
Therefore $D \in \mathfrak{sp}(6,\mathbb{R})$. Then, the $\mathrm{SU}(3)$-structure defined by
\begin{equation*}
\begin{aligned}
\omega&=e^{12}+e^{34}+e^{56},\\
\psi_+&=e^{135}-e^{146}-e^{236}-e^{245}.
\end{aligned}
\end{equation*}
is half-flat.
\end{example}

As a consequence of Propositions \ref{G2cocerrada} and \ref{G2cocerradainversa} we have the following:

\begin{theorem}
Let $\frh$ be a 6-dimensional Lie algebra and let $\frg$ be a 7-dimensional Lie algebra satisfying 
\begin{equation*}
\frg=\frh\oplus_D\mathbb{R}e_7,
\end{equation*}
with $D$ a derivation of $\mathfrak{h}$ given by \eqref{sp(6,R)} in terms of a basis $\{e_1,\dots,e_6\}$ of $\frh$. Then the following two conditions are equivalent:
\begin{enumerate}
\item The $\mathrm{SU}(3)$-structure on $\frh$ given by 
\begin{equation*}
\begin{aligned}
\omega&=e^{12}+e^{34}+e^{56},\\
\psi_+&=e^{135}-e^{146}-e^{236}-e^{245}.
\end{aligned}
\end{equation*}
is half-flat.
\item The $\mathrm{G}_2$-structure on $\frg$ given by 
\begin{equation*}
\varphi=e^{127}+e^{347}+e^{567}+e^{136}+e^{145}+e^{235}-e^{246},
\end{equation*}
is coclosed.
\end{enumerate}
\end{theorem}

\end{section}


\begin{section}{Construction of Lie algebras with closed $\mathrm{G}_2$-structures}
Finally, using the results of Section \ref{Closed} we describe all the 7-dimensional Lie algebras with closed $\mathrm{G}_2$-structures which are constructed as
\begin{equation*}
\frg=\frh\oplus_D\mathbb{R}e_7,
\end{equation*}
where $\frh$ denotes a 6-dimensional solvable Lie algebra with a symplectic half-flat structure and $D$ is a derivation of $\frh$. From \cite{FMOU} the 6-dimensional Lie algebras with a symplectic half-flat $\mathrm{SU}(3)$-structure are:
\begin{equation*}
\begin{aligned}
\fra&=(0,0,0,0,0,0);\\
\mathfrak{e}(1,1)\oplus \mathfrak{e}(1,1)&=(0,0,-e^{14},-e^{13},e^{25},-e^{26});\\
\frg_{5,1}\oplus \mathbb{R}&=(0,0,0,e^{15},0,e^{13});\\
\frg_{5,7^{-1,-1,1}}\oplus \mathbb{R}&=(-e^{15},e^{25},-e^{35},e^{45},0,0);\\
\frg_{5,17}^{\alpha, -\alpha, 1}\oplus \mathbb{R}&=(\alpha e^{15}+e^{35},-\alpha e^{25}+e^{45},-e^{15}+\alpha e^{35},-e^{25}-\alpha e^{45},0,0);\\
\frg_{6,N3}&=(0,e^{35},0,2e^{15},0,e^{13});\\
\frg_{6,38}^0&=(2e^{36},0,-e^{26},-e^{26}+e^{25},-e^{23}-e^{24},e^{23});\\
\frg_{6,54}^{0,-1}&=(e^{16}+e^{45},-e^{26},-e^{36}+e^{25},e^{46},0,0);\\
\frg_{6,118}^{0,-1,-1}&=(-e^{15}+e^{36},e^{46}+e^{25},-e^{16}-e^{35},-e^{45}-e^{26},0,0);\\
A_{6,13}^{-\frac{2}{3},\frac{1}{3},-1}&=(-\frac{1}{4}e^{14}-e^{23},\frac{1}{4}e^{24},-e^{26},-e^{26}+e^{25},-e^{23}-e^{24},e^{23});\\
A_{6,54}^{2,1}&=(-\frac{1}{2}e^{15},\frac{1}{2}e^{25}+e^{16},-\frac{1}{2}e^{35},\frac{1}{2}e^{45}+e^{36},0,-e^{56});\\
A_{6,70}^{\alpha,\frac{\alpha}{2}}(\alpha \neq 0)&=(-\frac{1}{2}e^{15}+\frac{1}{\alpha}e^{35}+e^{26}, \frac{1}{2}e^{25}+\frac{1}{\alpha}e^{45},-\frac{1}{\alpha}e^{15}-\frac{1}{2}e^{35}+e^{46},\\
& \quad  -\frac{1}{\alpha}e^{25}+\frac{1}{2}e^{45},0,e^{56});\\
A_{6,71}^{-\frac{3}{2}}&=(-\frac{3}{4}e^{16},\frac{3}{4}e^{26}+e^{35},\frac{1}{4}e^{36}+e^{45},-\frac{1}{4}e^{46}+e^{15},\frac{1}{2}e^{56},0);\\
N_{6,3}^{0,-2,0,2}&=(-2\cdot 3^{1/3}e^{16},2\cdot 3^{-1/3}e^{26},3^{-1/3}e^{36}+3^{2/3}e^{45},0,\\
& \quad -3^{-1/3}e^{34}-3^{-1/3}e^{56},0).\\
\end{aligned}
\end{equation*}
The structure equations of the previously mentioned Lie algebras are given in terms of an adapted basis, that is, a basis such that the forms 
\begin{equation*}
\begin{aligned}
\omega&=e^{12}+e^{34}+e^{56},\\
\psi_+&=e^{135}-e^{146}-e^{236}-e^{245},
\end{aligned}
\end{equation*}
are closed and therefore describe a symplectic half-flat $\mathrm{SU}(3)$-structure. 

\begin{proposition}\label{proposition}
The Lie algebras described in Table 1 of the Appendix admit the closed $\mathrm{G}_2$-structure given by the 3-form \eqref{G2form}, that is
\begin{equation*}
\varphi=e^{127}+e^{347}+e^{567}+e^{135}-e^{146}-e^{236}-e^{245}.
\end{equation*}
\begin{proof}
For each one of the 6-dimensional solvable Lie algebras admitting a symplectic half-flat $\mathrm{SU}(3)$-structure, namely $\frh$ we consider the Lie algebras
\begin{equation*}
\frg=\frh\oplus_D \mathbb{R}e_7,
\end{equation*}
with $D$ being the real representation of certain $A \in \mathfrak{sl}(3,\mathbb{C})$, that is, $D$ is of the form 
\begin{equation*}
 D = \left(
\begin{array}{cc|cc|cc}
 a_{1,1} & a_{1,2} & a_{1,3} & a_{1,4} & a_{1,5} & a_{1,6}  \\
 -a_{1,2} & a_{1,1} & -a_{1,4} & a_{1,3} & -a_{1,6} & a_{1,5}  \\\hline
 a_{3,1} & a_{3,2} & a_{3,3} & a_{3,4} & a_{3,5} & a_{3,6}  \\
 -a_{3,2} & a_{3,1} & -a_{3,4} & a_{3,3} & -a_{3,6} & a_{3,5}  \\\hline
 a_{5,1} & a_{5,2} & a_{5,3} & a_{5,4} & -a_{1,1}-a_{3,3} & -a_{1,2}-a_{3,4}  \\
 -a_{5,2} & a_{5,1} & -a_{5,4} & a_{5,3} & a_{1,2}+a_{3,4} & -a_{1,1}-a_{3,3}  \\
\end{array}
\right).
\end{equation*}
Then, we compute for which values of the parameters $a_{i,j}$, the matrix $D$ represents a derivation of $\mathfrak{h}$ and thus, from Proposition 1.1 the 3-form
\begin{equation*}
\varphi= e^{127}+e^{347}+e^{567}+e^{135}-e^{146}-e^{236}-e^{245}
\end{equation*}
defines a closed $\mathrm{G}_2$-structure on $\mathfrak{g}=\mathfrak{h}\oplus_{D} \mathbb{R}e_7$.

\medskip

\begin{itemize}
\item \textbf{ $\frh=\mathfrak{a}$}
\end{itemize}
The structure equations of 
\begin{equation*}
\mathfrak{g}=\mathfrak{a}\oplus_{D} \mathbb{R}e_7,
\end{equation*}
with $D$ a derivation described as in \eqref{sl(3,C)} are
\begin{equation}\label{abeliana}
\begin{aligned}
de^1= &\, a_{1,1} e^{17}-a_{1,2} e^{27}+a_{3,1} e^{37}-a_{3,2} e^{47}+a_{5,1} e^{57}-a_{5,2} e^{67},\\
de^2= &\, a_{1,2} e^{17}+a_{1,1} e^{27}+a_{3,2} e^{37}+a_{3,1} e^{47}+a_{5,2} e^{57}+a_{5,1} e^{67},\\
de^3= &\, a_{1,3} e^{17}-a_{1,4} e^{27}+a_{3,3} e^{37}-a_{3,4} e^{47}+a_{5,3} e^{57}-a_{5,4} e^{67},\\
de^4= &\, a_{1,4} e^{17}+a_{1,3} e^{27}+a_{3,4} e^{37}+a_{3,3} e^{47}+a_{5,4} e^{57}+a_{5,3} e^{67},\\
de^5= &\, a_{1,5} e^{17}-a_{1,6} e^{27}+a_{3,5} e^{37}-a_{3,6} e^{47} +\left(-a_{1,1}-a_{3,3}\right) e^{57}+\left(a_{1,2}+a_{3,4}\right) e^{67},\\
de^6= &\, a_{1,6} e^{17}+a_{1,5} e^{27}+a_{3,6} e^{37}+a_{3,5} e^{47}+\left(-a_{1,2}-a_{3,4}\right) e^{57}+\left(-a_{1,1}-a_{3,3}\right) e^{67},\\
de^7= &\, 0.
\end{aligned}
\end{equation}
The condition of $D$ being a derivation of $\fra$ is equivalent to the vanishing of the differential operator when applied twice. From \eqref{abeliana}   
\begin{align*}
d^2e^1= &\, 0, \qquad d^2e^2= \, 0, \qquad d^2e^3= \, 0, \qquad d^2e^4= \, 0, \\
d^2e^5= &\, 0, \qquad d^2e^6= \, 0, \qquad d^2e^7= \, 0, \\
\end{align*}
and therefore $D$ is derivation of $\mathfrak{a}$ for every $a_{1,1}, a_{1,2}, a_{1,3}, a_{1,4}, a_{1,5}, a_{1,6}, a_{3,1}, a_{3,2}, $ $a_{3,3}, a_{3,4}, a_{3,5}, a_{3,6}, a_{5,1}, a_{5,2}, a_{5,3}, a_{5,4} \in \mathbb{R}$. Thus, the Lie algebra $\frg$ which structure equations are described in \eqref{abeliana} admits the closed $\mathrm{G}_2$-structure $\eqref{G2form}$.

\medskip

\begin{itemize}
\item \textbf{ $\frh=\mathfrak{e}(1,1)\oplus\mathfrak{e}(1,1)$}
\end{itemize}
The structure equations of 
\begin{equation*}
\mathfrak{g}=\big(\mathfrak{e}(1,1)\oplus\mathfrak{e}(1,1)\big)\oplus_{D} \mathbb{R}e_7,
\end{equation*}
with $D$ a derivation described as in \eqref{sl(3,C)} are
\begin{equation*}
\begin{aligned}
de^1= &\, a_{1,1} e^{17}-a_{1,2} e^{27}+a_{3,1} e^{37}-a_{3,2} e^{47}+a_{5,1} e^{57}-a_{5,2} e^{67},\\
de^2= &\, a_{1,2} e^{17}+a_{1,1} e^{27}+a_{3,2} e^{37}+a_{3,1} e^{47}+a_{5,2} e^{57}+a_{5,1} e^{67},\\
de^3= &\, -e^{14}+a_{1,3} e^{17}-a_{1,4} e^{27}+a_{3,3} e^{37}-a_{3,4} e^{47}+a_{5,3} e^{57}-a_{5,4} e^{67},\\
de^4= &\, -e^{13}+a_{1,4} e^{17}+a_{1,3} e^{27}+a_{3,4} e^{37}+a_{3,3} e^{47}+a_{5,4} e^{57}+a_{5,3} e^{67},\\
de^5= &\, e^{25}+a_{1,5} e^{17}-a_{1,6} e^{27}+a_{3,5} e^{37}-a_{3,6} e^{47}\\
           &\,+\left(-a_{1,1}-a_{3,3}\right) e^{57}+\left(a_{1,2}+a_{3,4}\right) e^{67},\\
de^6= &\, -e^{26}+a_{1,6} e^{17}+a_{1,5} e^{27}+a_{3,6} e^{37}+a_{3,5} e^{47}\\
           &\,+\left(-a_{1,2}-a_{3,4}\right) e^{57}+\left(-a_{1,1}-a_{3,3}\right) e^{67},\\
de^7= &\, 0.
\end{aligned}
\end{equation*}
The condition of $D$ being a derivation of $\mathfrak{e}(1,1)\oplus\mathfrak{e}(1,1)$ is equivalent to the vanishing of the differential operator when applied twice. Thus 
\begin{equation*}
\begin{aligned}
d^2e^1= &\, a_{3,2} e^{137}-a_{3,1} e^{147}+a_{5,1} e^{257}+a_{5,2} e^{267},\\
d^2e^2= &\, -a_{3,1} e^{137}-a_{3,2} e^{147}+a_{5,2} e^{257}-a_{5,1} e^{267},\\
d^2e^3= &\, a_{1,3} e^{127}+2 a_{3,4} e^{137}+a_{1,1} e^{147}+a_{5,4} e^{157} +a_{5,3} e^{167}-a_{1,2} e^{247}+a_{5,3} e^{257}\\
               &\, +a_{5,4} e^{267}+a_{3,1} e^{347}-a_{5,1} e^{457}+a_{5,2} e^{467},\\
d^2e^4= &\, -a_{1,4} e^{127}+a_{1,1} e^{137}-2 a_{3,4} e^{147}+a_{5,3} e^{157}-a_{5,4} e^{167}-a_{1,2} e^{237}+a_{5,4} e^{257}\\
               &\, -a_{5,3} e^{267}+a_{3,2} e^{347}-a_{5,1} e^{357}+a_{5,2} e^{367},\\
d^2e^5= &\, a_{1,5} e^{127}+a_{3,6} e^{137}-a_{3,5} e^{147}-a_{1,2} e^{157}-a_{3,5} e^{237}+a_{3,6} e^{247}\\
               &\, -a_{1,1} e^{257} +\left(-2 a_{1,2}-2 a_{3,4}\right) e^{267}-a_{3,2} e^{357}-a_{3,1} e^{457}+a_{5,1} e^{567},\\
d^2e^6= &\, -a_{1,6} e^{127}-a_{3,5} e^{137}-a_{3,6} e^{147}+a_{1,2} e^{167}+a_{3,6} e^{237}+a_{3,5} e^{247}\\
               &\, +\left(-2 a_{1,2}-2 a_{3,4}\right) e^{257}+a_{1,1} e^{267}+a_{3,2} e^{367}+a_{3,1} e^{467}+a_{5,2} e^{567},\\
d^2e^7= &\, 0,
\end{aligned}
\end{equation*}
and after solving the system $d^2=0$ we conclude that the derivations $D \in \mathfrak{sl}(3, \mathbb{C})$ of $\mathfrak{e}(1,1)\oplus\mathfrak{e}(1,1)$ are of the form 
\begin{equation*}
 D = \left(
\begin{array}{cc|cc|cc}
&  & &&&  \\
 & & &&& \\ \hline
  & & a_{33} &  &&  \\
 & &  & a_{33} & & \\ \hline
 & &&& -a_{33} &   \\
 &&&&& -a_{33}  \\
\end{array}
\right).
\end{equation*}
Therefore, the family of Lie algebras $\frg=\big(\mathfrak{e}(1,1)\oplus\mathfrak{e}(1,1))\oplus_{D} \mathbb{R}e_7$ whose structure equations are 
\begin{align*}
\big(\mathfrak{e}(1,1)\oplus\mathfrak{e}(1,1))\oplus_{D} \mathbb{R}e_7 = & (0,0,-e^{14}+a_{3,3} e^{37}, -e^{13}+a_{3,3} e^{47}, e^{25}-a_{3,3} e^{57}, \\ &-e^{26}-a_{3,3} e^{67})   
\end{align*}
admits the closed $\mathrm{G}_2$-structure defined by \eqref{G2form}.

\medskip


\begin{itemize}
\item \textbf{ $\frh=\mathfrak{g}_{5,1}\oplus\mathbb{R}$}
\end{itemize}
The structure equations of 
\begin{equation*}
\mathfrak{g}=\big(\mathfrak{g}_{5,1}\oplus\mathbb{R}\big)\oplus_{D} \mathbb{R}e_7,
\end{equation*}
with $D$ a derivation described as in \eqref{sl(3,C)} are
\begin{equation*}
\begin{aligned}
de^1= &\, a_{1,1} e^{17}-a_{1,2} e^{27}+a_{3,1} e^{37}-a_{3,2} e^{47}+a_{5,1} e^{57}-a_{5,2} e^{67},\\
de^2= &\, a_{1,2} e^{17}+a_{1,1} e^{27}+a_{3,2} e^{37}+a_{3,1} e^{47}+a_{5,2} e^{57}+a_{5,1} e^{67},\\
de^3= &\, a_{1,3} e^{17}-a_{1,4} e^{27}+a_{3,3} e^{37}-a_{3,4} e^{47}+a_{5,3} e^{57}-a_{5,4} e^{67},\\
de^4= &\, e^{15}+a_{1,4} e^{17}+a_{1,3} e^{27}+a_{3,4} e^{37}+a_{3,3} e^{47}+a_{5,4} e^{57}+a_{5,3} e^{67},\\
de^5= &\, a_{1,5} e^{17}-a_{1,6} e^{27}+a_{3,5} e^{37}-a_{3,6} e^{47}+\left(-a_{1,1}-a_{3,3}\right) e^{57}+\left(a_{1,2}+a_{3,4}\right) e^{67},\\
de^6= &\, e^{13}+a_{1,6} e^{17}+a_{1,5} e^{27}+a_{3,6} e^{37}+a_{3,5} e^{47}+\left(-a_{1,2}-a_{3,4}\right) e^{57}+\left(-a_{1,1}-a_{3,3}\right) e^{67},\\
de^7= &\, 0.
\end{aligned}
\end{equation*}
The condition of $D$ being a derivation of $\mathfrak{g}_{5,1}\oplus\mathbb{R}$ is equivalent to the vanishing of the differential operator when applied twice. Thus 
\begin{equation*}
\begin{aligned}
d^2e^1= &\,-a_{5,2} e^{137}-a_{3,2} e^{157},\\
d^2e^2= &\, a_{5,1} e^{137}+a_{3,1} e^{157},\\
d^2e^3= &\, -a_{5,4} e^{137}-a_{3,4} e^{157},\\
d^2e^4= &\, a_{1,6} e^{127}+\left(a_{5,3}-a_{3,5}\right) e^{137}+a_{3,6} e^{147}+2 a_{3,3} e^{157}+\left(-a_{1,2}-a_{3,4}\right) e^{167}\\
               &\, +a_{1,2} e^{257}-a_{3,1} e^{357}+a_{3,2} e^{457}-a_{5,2} e^{567}\\
d^2e^5= &\, \left(a_{1,2}+a_{3,4}\right) e^{137}-a_{3,6} e^{157},\\
d^2e^6= &\, a_{1,4} e^{127}+\left(-2 a_{1,1}-2 a_{3,3}\right) e^{137}+a_{3,4} e^{147}+\left(a_{3,5}-a_{5,3}\right) e^{157}+a_{5,4} e^{167}\\
               &\,+a_{1,2} e^{237}-a_{3,2} e^{347}+a_{5,1} e^{357}-a_{5,2} e^{367},\\
d^2e^7= &\, 0,
\end{aligned}
\end{equation*}
and after solving the system $d^2=0$ we conclude that the derivations $D \in \mathfrak{sl}(3, \mathbb{C})$ of $\mathfrak{g}_{5,1}\oplus\mathbb{R}$ are of the form 
\begin{equation*}
 D = \left(
\begin{array}{cc|cc|cc}
&& a_{1,3} && a_{1,5} & \\
&&& a_{1,3} && a_{1,5}  \\ \hline
&&&& a_{3,5} &  \\
&&&&& a_{3,5} \\ \hline
 & &a_{3,5}&&  &   \\
 &&&a_{3,5}&& \\
\end{array}
\right).
\end{equation*}
\smallskip
Therefore, the family of Lie algebras $\frg=\big(\mathfrak{g}_{5,1}\oplus\mathbb{R}\big)\oplus_{D} \mathbb{R}e_7$ whose structure equations are 
\begin{align*}
\big(\mathfrak{g}_{5,1}\oplus\mathbb{R}\big)\oplus_{D} \mathbb{R}e_7 = & (0,0,a_{1,3}e^{17}+a_{3,5} e^{57}, e^{15}+a_{1,3} e^{27}+a_{3,5}e^{67}, a_{1,5} e^{17}+a_{3,5}e^{37}, \\
& e^{13}+a_{1,5} e^{27}+a_{3,5}e^{47}, 0 )   
\end{align*}
are such that the $\mathrm{G}_2$-form \eqref{G2form} is closed.

\medskip


\begin{itemize}
\item \textbf{ $\frh=\mathfrak{g}_{5,7}^{-1,-1,1}\oplus\mathbb{R}$}
\end{itemize}
The structure equations of 
\begin{equation*}
\mathfrak{g}=\big(\mathfrak{g}_{5,7}^{-1,-1,1}\oplus\mathbb{R}\big)\oplus_{D} \mathbb{R}e_7,
\end{equation*}
with $D$ a derivation described as in \eqref{sl(3,C)} are
\begin{equation*}
\begin{aligned}
de^1= &\, -e^{15}+a_{1,1} e^{17}-a_{1,2} e^{27}+a_{3,1} e^{37}-a_{3,2} e^{47}+a_{5,1} e^{57}-a_{5,2} e^{67},\\
de^2= &\, e^{25}+a_{1,2} e^{17}+a_{1,1} e^{27}+a_{3,2} e^{37}+a_{3,1} e^{47}+a_{5,2} e^{57}+a_{5,1} e^{67},\\
de^3= &\, -e^{35}+a_{1,3} e^{17}-a_{1,4} e^{27}+a_{3,3} e^{37}-a_{3,4} e^{47}+a_{5,3} e^{57}-a_{5,4} e^{67},\\
de^4= &\, e^{45}+a_{1,4} e^{17}+a_{1,3} e^{27}+a_{3,4} e^{37}+a_{3,3} e^{47}+a_{5,4} e^{57}+a_{5,3} e^{67},\\
de^5= &\,  a_{1,5}e^{17}-a_{1,6} e^{27}+a_{3,5} e^{37}-a_{3,6} e^{47}+\left(-a_{1,1}-a_{3,3}\right) e^{57}+\left(a_{1,2}+a_{3,4}\right) e^{67},\\
de^6= &\, a_{1,6} e^{17}+a_{1,5} e^{27}+a_{3,6} e^{37}+a_{3,5} e^{47}+\left(-a_{1,2}-a_{3,4}\right) e^{57}+\left(-a_{1,1}-a_{3,3}\right) e^{67},\\
de^7= &\, 0.
\end{aligned}
\end{equation*}
Thus after imposing $d^2=0$ we conclude that the real representations of the derivations $D \in \mathfrak{sl}(3, \mathbb{C})$ of $\mathfrak{g}_{5,7}^{-1,-1,1}\oplus\mathbb{R}$ are 
\begin{equation*}
 D = \left(
\begin{array}{cc|cc|cc}
 a_{1,1} && a_{1,3} &&& \\
& a_{1,1} && a_{1,3} && \\ \hline
a_{3,1} && -a_{1,1} &&&  \\
& a_{3,1} && -a_{1,1} && \\ \hline
 & &&&  &   \\
 &&&&& \\
\end{array}
\right).
\end{equation*}
Hence, the family of Lie algebras $\frg=\big(\mathfrak{g}_{5,7}^{-1,-1,1}\oplus\mathbb{R}\big)\oplus_{D} \mathbb{R}e_7$ whose structure equations are 
\begin{align*}
\big(\mathfrak{g}_{5,7}^{-1,-1,1}\oplus\mathbb{R}\big)\oplus_{D} \mathbb{R}e_7 = & (-e^{15}+a_{1,3}e^{17}-a_{1,1} e^{37}, e^{25}+a_{1,1} e^{17}+a_{3,1}e^{37}, \\
 & \,-e^{35}+a_{1,3} e^{17}-a_{1,1}e^{37},  e^{45}+a_{1,3} e^{27}-a_{1,1}e^{47}, 0, 0, 0)   
\end{align*}
is such that the $\mathrm{G}_2$-form \eqref{G2form} is closed for all $a_{1,1}, a_{1,3}$ and $a_{3,1}$ real numbers.

\medskip


\begin{itemize}
\item \textbf{ $\frh=\mathfrak{g}_{5,17}^{\alpha,-\alpha,1}\oplus\mathbb{R}$ with $ \alpha \geq 0$}
\end{itemize}
The structure equations of 
\begin{equation*}
\mathfrak{g}=\big(\mathfrak{g}_{5,17}^{\alpha,-\alpha,1}\oplus\mathbb{R}\big)\oplus_{D} \mathbb{R}e_7,
\end{equation*}
with $D$ as in \eqref{sl(3,C)} are
\begin{equation*}
\begin{aligned}
de^1= &\, \alpha e^{15}+ e^{35}+a_{1,1} e^{17}-a_{1,2} e^{27}+a_{3,1} e^{37}-a_{3,2} e^{47}+a_{5,1} e^{57}-a_{5,2} e^{67},\\
de^2= &\, -\alpha e^{25}+e^{45}+a_{1,2} e^{17}+a_{1,1} e^{27}+a_{3,2} e^{37}+a_{3,1} e^{47}+a_{5,2} e^{57}+a_{5,1} e^{67},\\
de^3= &\, -e^{15}+\alpha e^{35}+a_{1,3} e^{17}-a_{1,4} e^{27}+a_{3,3} e^{37}-a_{3,4} e^{47}+a_{5,3} e^{57}-a_{5,4} e^{67},\\
de^4= &\, -e^{25}-\alpha e^{45}+a_{1,4} e^{17}+a_{1,3} e^{27}+a_{3,4} e^{37}+a_{3,3} e^{47}+a_{5,4} e^{57}+a_{5,3} e^{67},\\
de^5= &\,  a_{1,5}e^{17}-a_{1,6} e^{27}+a_{3,5} e^{37}-a_{3,6} e^{47}+\left(-a_{1,1}-a_{3,3}\right) e^{57}+\left(a_{1,2}+a_{3,4}\right) e^{67},\\
de^6= &\, a_{1,6} e^{17}+a_{1,5} e^{27}+a_{3,6} e^{37}+a_{3,5} e^{47}+\left(-a_{1,2}-a_{3,4}\right) e^{57}+\left(-a_{1,1}-a_{3,3}\right) e^{67},\\
de^7= &\, 0.
\end{aligned}
\end{equation*}
As before we impose the condition $d^2=0$, obtaining that the derivations $D \in \mathfrak{sl}(3, \mathbb{C})$ of $\mathfrak{g}_{5,17}^{\alpha,-\alpha,1}\oplus\mathbb{R}$ are 
\begin{equation*}
 D = \left(
\begin{array}{cc|cc|cc}
 && a_{1,3} &&& \\
&  && a_{1,3} && \\ \hline
-a_{1,3}&&  &&&  \\
& -a_{1,3}&&  && \\ \hline
 & &&&  &   \\
 &&&&& \\
\end{array}
\right).
\end{equation*}
Therefore, the family of Lie algebras $\frg=\big(\mathfrak{g}_{5,17}^{\alpha,-\alpha,1}\oplus\mathbb{R}\big)\oplus_{D} \mathbb{R}e_7$ with $\alpha \geq 0$, whose structure equations are 
\begin{align*}
\big(\mathfrak{g}_{5,17}^{\alpha,-\alpha,1}\oplus\mathbb{R}\big)\oplus_{D} \mathbb{R}e_7 = & (\alpha e^{15}+e^{35}-a_{1,3}e^{37}, -\alpha e^{25}+ e^{45}-a_{1,3}e^{47}, \\
 & \,-e^{15}+\alpha e^{35}+a_{1,3}e^{17},  -e^{25}- \alpha e^{45}+a_{1,3}e^{27}, 0, 0, 0)   
\end{align*}
satisfies that the $\mathrm{G}_2$-form described by \eqref{G2form} is closed for all $\alpha \geq 0$ and $a_{1,3}$ real number.

\medskip


\begin{itemize}
\item \textbf{ $\frh=\mathfrak{g}_{6,N3}$}
\end{itemize}
The structure equations of 
\begin{equation*}
\mathfrak{g}=\mathfrak{g}_{6,N3}\oplus_{D} \mathbb{R}e_7,
\end{equation*}
with $D$ the real representation of certain $A \in \mathfrak{sl}(3,\mathbb{C})$ are
\begin{equation*}
\begin{aligned}
de^1= &\, a_{1,1} e^{17}-a_{1,2} e^{27}+a_{3,1} e^{37}-a_{3,2} e^{47}+a_{5,1} e^{57}-a_{5,2} e^{67},\\
de^2= &\, e^{35}+a_{1,2} e^{17}+a_{1,1} e^{27}+a_{3,2} e^{37}+a_{3,1} e^{47}+a_{5,2} e^{57}+a_{5,1} e^{67},\\
de^3= &\, a_{1,3} e^{17}-a_{1,4} e^{27}+a_{3,3} e^{37}-a_{3,4} e^{47}+a_{5,3} e^{57}-a_{5,4} e^{67},\\
de^4= &\, 2e^{15}+a_{1,4} e^{17}+a_{1,3} e^{27}+a_{3,4} e^{37}+a_{3,3} e^{47}+a_{5,4} e^{57}+a_{5,3} e^{67},\\
de^5= &\, a_{1,5} e^{17}-a_{1,6} e^{27}+a_{3,5} e^{37}-a_{3,6} e^{47}+\left(-a_{1,1}-a_{3,3}\right) e^{57}+\left(a_{1,2}+a_{3,4}\right) e^{67},\\
de^6= &\, e^{13}+a_{1,6} e^{17}+a_{1,5} e^{27}+a_{3,6} e^{37}+a_{3,5} e^{47}+\left(-a_{1,2}-a_{3,4}\right) e^{57}\\
           &\, +\left(-a_{1,1}-a_{3,3}\right) e^{67},\\
de^7= &\, 0.\end{aligned}
\end{equation*}
Proceeding as before we obtain that the derivation $D$ is described by\begin{equation*}
 D = \left(
\begin{array}{cc|cc|cc}
 && a_{1,3} && a_{1,5}& \\
&  && a_{1,3} && a_{1,5} \\ \hline
\frac{a_{1,3}}{2}&&  && a_{3,5}&  \\
& \frac{a_{1,3}}{2}&&  && a_{3,5} \\ \hline
 -a_{1,5} && 2a_{3,5}&&  &   \\
 & -a_{1,5} && 2a_{3,5}&& \\
\end{array}
\right).
\end{equation*}
Thus, the family of Lie algebras $\frg=\mathfrak{g}_{6,N3}\oplus_{D} \mathbb{R}e_7$ has structure equations
\begin{align*}
\mathfrak{g}_{6,N3}\oplus_{D} \mathbb{R}e_7 = & \big(\frac{a_{1,3}}{2}e^{37}-a_{1,5}e^{57}, e^{35}+\frac{a_{1,3}}{2} e^{47}-a_{1,5}e^{67},  a_{1,3}e^{17}+2a_{3,5}e^{57}, \\
 & \, a_{1,3}e^{27}+2a_{3,5}e^{67}, a_{1,5}e^{17}+a_{3,5}e^{37},  e^{13}+a_{1,5} e^{27}+a_{3,5}e^{47}, 0\big).
\end{align*}
Hence, the $\mathrm{G}_2$ form \eqref{G2form} is closed for all $a_{1,3}, a_{1,5}$ and $a_{3,5}$ real numbers.

\medskip


\begin{itemize}
\item \textbf{ $\frg=\mathfrak{g}_{6,38}^0$}
\end{itemize}
The structure equations of 
\begin{equation*}
\mathfrak{g}=\mathfrak{g}_{6,38}^0\oplus_{D} \mathbb{R}e_7,
\end{equation*}
with $D$ as in \eqref{sl(3,C)} are
\begin{equation*}
\begin{aligned}
de^1= &\, 2e^{36}+a_{1,1} e^{17}-a_{1,2} e^{27}+a_{3,1} e^{37}-a_{3,2} e^{47}+a_{5,1} e^{57}-a_{5,2} e^{67},\\
de^2= &\, a_{1,2} e^{17}+a_{1,1} e^{27}+a_{3,2} e^{37}+a_{3,1} e^{47}+a_{5,2} e^{57}+a_{5,1} e^{67},\\
de^3= &\, -e^{26}+a_{1,3} e^{17}-a_{1,4} e^{27}+a_{3,3} e^{37}-a_{3,4} e^{47}+a_{5,3} e^{57}-a_{5,4} e^{67},\\
de^4= &\, -e^{26}+e^{25}+a_{1,4} e^{17}+a_{1,3} e^{27}+a_{3,4} e^{37}+a_{3,3} e^{47}+a_{5,4} e^{57}+a_{5,3} e^{67},\\
de^5= &\, -e^{23}-e^{24}+a_{1,5} e^{17}-a_{1,6} e^{27}+a_{3,5} e^{37}-a_{3,6} e^{47}+\left(-a_{1,1}-a_{3,3}\right) e^{57}\\
           &\, +\left(a_{1,2}+a_{3,4}\right) e^{67},\\
de^6= &\, e^{23}+a_{1,6} e^{17}+a_{1,5} e^{27}+a_{3,6} e^{37}+a_{3,5} e^{47}+\left(-a_{1,2}-a_{3,4}\right) e^{57}\\
           &\, +\left(-a_{1,1}-a_{3,3}\right) e^{67},\\
de^7= &\, 0.\end{aligned}
\end{equation*}
Impossing $d^2=0$ we obtain that the unique derivation $D$ of $\mathfrak{g}_{6,38}^0$ of the form \eqref{sl(3,C)} is given by the null matrix. Thus, $\mathfrak{g}$ is direct sum, that is

\begin{align*}
\mathfrak{g}_{6,38}^0\oplus \mathbb{R}e_7 = & (2e^{36}, 0, -e^{26}, -e^{26}+e^{25}, -e^{23}-e^{24}, e^{23},0).
\end{align*}
It is clear that for this Lie algebra the $\mathrm{G}_2$ form \eqref{G2form} is closed.

\medskip


\begin{itemize}
\item \textbf{ $\frh=\mathfrak{g}_{6,54}^{0,-1}$}
\end{itemize}
Exactly as before we have that 
\begin{equation*}
\mathfrak{g}=\mathfrak{g}_{6,54}^{0,-1}\oplus_{D} \mathbb{R}e_7,
\end{equation*}
with $D$ as in \eqref{sl(3,C)} is exactly a direct sum, that is
\begin{align*}
\mathfrak{g}_{6,54}^{0,-1}\oplus \mathbb{R}e_7 = & (e^{16}+e^{45}, -e^{26}, -e^{36}+e^{25}, e^{46}, 0, 0,0),
\end{align*}
and therefore the $\mathrm{G}_2$ form \eqref{G2form} is closed for this algebra.

\medskip


\begin{itemize}
\item \textbf{ $\frh=\mathfrak{g}_{6,118}^{0,-1,-1}$}
\end{itemize}
The structure equations of 
\begin{equation*}
\mathfrak{g}=\mathfrak{g}_{6,118}^{0,-1,-1}\oplus_{D} \mathbb{R}e_7,
\end{equation*}
with $D$ the real representation of certain $3 \times 3$ complex matrix without trace are
\begin{equation*}
\begin{aligned}
de^1= &\, -e^{15}+e^{36}+a_{1,1} e^{17}-a_{1,2} e^{27}+a_{3,1} e^{37}-a_{3,2} e^{47}+a_{5,1} e^{57}-a_{5,2} e^{67},\\
de^2= &\, e^{46}+e^{25}+a_{1,2} e^{17}+a_{1,1} e^{27}+a_{3,2} e^{37}+a_{3,1} e^{47}+a_{5,2} e^{57}+a_{5,1} e^{67},\\
de^3= &\, -e^{16}-e^{35}+a_{1,3} e^{17}-a_{1,4} e^{27}+a_{3,3} e^{37}-a_{3,4} e^{47}+a_{5,3} e^{57}-a_{5,4} e^{67},\\
de^4= &\, e^{45}-e^{26}+a_{1,4} e^{17}+a_{1,3} e^{27}+a_{3,4} e^{37}+a_{3,3} e^{47}+a_{5,4} e^{57}+a_{5,3} e^{67},\\
de^5= &\,a_{1,5} e^{17}-a_{1,6} e^{27}+a_{3,5} e^{37}-a_{3,6} e^{47}+\left(-a_{1,1}-a_{3,3}\right) e^{57}\\
           &\, +\left(a_{1,2}+a_{3,4}\right) e^{67},\\
de^6= &\,a_{1,6} e^{17}+a_{1,5} e^{27}+a_{3,6} e^{37}+a_{3,5} e^{47}+\left(-a_{1,2}-a_{3,4}\right) e^{57}\\
           &\, +\left(-a_{1,1}-a_{3,3}\right) e^{67},\\
de^7= &\, 0.\end{aligned}
\end{equation*}
Solving the equation obtained from the fact that $D$ has to be a derivation of $\mathfrak{g}_{6,118}^{0,-1,-1}$ we have that $D \in \mathfrak{sl}(3,\mathbb{C})$ is such that
\begin{equation*}
 D = \left(
\begin{array}{cc|cc|cc}
 && a_{1,3} && & \\
&  && a_{1,3} && \\ \hline
-a_{1,3}&&  && &  \\
& -a_{1,3}&&  &&  \\ \hline
&&&&  &   \\
 & &&&& \\
\end{array}
\right).
\end{equation*}
Thus, from Proposition 1.1 the family of Lie algebras with structure equations 
\begin{align*}
\mathfrak{g}_{6,118}^{0,-1,-1}\oplus_D \mathbb{R}e_7 = & (-e^{15}+e^{36}-a_{1,3}e^{37}, e^{46}+e^{25}-a_{1,3}e^{47}, -e^{16}-e^{35}+a_{1,3}e^{17},\\
& \, e^{45}-e^{26}-a_{1,3}e^{27}, 0, 0, 0)
\end{align*}
is such that the $\mathrm{G}_2$ form \eqref{G2form} is closed for all $a_{1,3}$ real number.

\medskip


\begin{itemize}
\item \textbf{ $\frh=A_{6,13}^{-\frac{2}{3},\frac{1}{3},-1}$}
\end{itemize}
With the same procedure as for the previous Lie algebras we have that a derivation $D$ of the Lie algebra $A_{6,13}^{-\frac{2}{3},\frac{1}{3},-1}$ being the real representation of a $3 \times 3$ complex matrix without trace has to be of the form
\begin{equation*}
 D = \left(
\begin{array}{cc|cc|cc}
a_{1,1} &&  && & \\
&  a_{1,1} && && \\ \hline
&&  && &  \\
& &&  &&  \\ \hline
&&&&  -a_{1,1} &   \\
 & &&&& -a_{1,1} \\
\end{array}
\right).
\end{equation*}
Therefore, from Proposition 1.1 the family of Lie algebras with structure equations 
\begin{align*}
A_{6,13}^{-\frac{2}{3},\frac{1}{3},-1}\oplus_D \mathbb{R}e_7 = & \big(-\frac{1}{4}e^{14}-e^{23}+a_{1,1}e^{17}, \frac{1}{4}e^{24}+a_{1,1}e^{27}, -\frac{1}{2}e^{34}, 0, \\ 
& \,-\frac{3}{4} e^{45}-a_{1,1}e^{57}, \frac{3}{4} e^{46}-a_{1,1}e^{67}, 0 \big)
\end{align*}
is such that the $\mathrm{G}_2$ form \eqref{G2form} is closed for all $a_{1,1}$ real number.

\medskip


\begin{itemize}
\item \textbf{ $\frh=A_{6,54}^{2,1}$}
\end{itemize}
For this Lie algebra can be checked that the derivations $D$ of $A_{6,54}^{2,1}$ of the form \eqref{sl(3,C)} are such that
\begin{equation*}
 D = \left(
\begin{array}{cc|cc|cc}
a_{1,1} && a_{1,3} && & \\
&  a_{1,1} && a_{1,3} && \\ \hline
a_{3,1} && -a_{1,1} && &  \\
& a_{3,1} && -a_{1,1} &&  \\ \hline
&&&&&   \\
 & &&&& \\
\end{array}
\right).
\end{equation*}
Henceforth, the family of Lie algebras with structure equations 
\begin{align*}
A_{6,54}^{2,1}\oplus_D \mathbb{R}e_7 = & \big(-\frac{1}{2}e^{15}+a_{1,1}e^{17}+a_{3,1}e^{37}, \frac{1}{2}e^{25}+e^{16}+a_{1,1}e^{27}+a_{3,1}e^{47}, \\
& \, -\frac{1}{2}e^{35}+a_{1,3}e^{17}-a_{1,1}e^{37}, -\frac{1}{2}e^{45}+e^{36}+a_{1,3}e^{27}-a_{1,1}e^{47}, 0, e^{56}, 0 \big)
\end{align*}
is such that the $\mathrm{G}_2$ form \eqref{G2form} is closed for all $a_{1,1}, a_{1,3}$ and $a_{3,1}$ real numbers.

\medskip


\begin{itemize}
\item \textbf{ $\frh=A_{6,70}^{\alpha,\frac{\alpha}{2}} \, (\alpha \neq 0)$}
\end{itemize}
For this Lie algebra the derivations  of the form \eqref{sl(3,C)} are such that
\begin{equation*}
 D = \left(
\begin{array}{cc|cc|cc}
 && a_{1,3} && & \\
&   && a_{1,3} && \\ \hline
-a_{1,3} && && &  \\
& -a_{1,3} && &&  \\ \hline
&&&&&   \\
 & &&&& \\
\end{array}
\right).
\end{equation*}
Therefore, the family of Lie algebras with structure equations 
\begin{align*}
A_{6,70}^{\alpha,\frac{\alpha}{2}}\oplus_D \mathbb{R}e_7 = & \big(-\frac{1}{2}e^{15}+\frac{1}{\alpha}e^{35}+e^{26}-a_{1,3}e^{37}, \frac{1}{2}e^{25}+\frac{1}{\alpha}e^{45}-a_{1,3}e^{47},\\
& \,
-\frac{1}{\alpha}e^{15}-\frac{1}{2}e^{35}+e^{46}+a_{1,3}e^{17}, -\frac{1}{\alpha}e^{25}+\frac{1}{2}e^{45}+a_{1,3}e^{27}, 0, e^{56}, 0 \big)
\end{align*}
is such that the $\mathrm{G}_2$ form \eqref{G2form} is closed for any $a_{1,3}\in \mathbb{R}$.  

\medskip


\begin{itemize}
\item \textbf{ $\frh=A_{6,71}^{-\frac{3}{2}}$}
\end{itemize}
The structure equations of 
\begin{equation*}
\mathfrak{g}=A_{6,71}^{-\frac{3}{2}}\oplus_{D} \mathbb{R}e_7,
\end{equation*}
with $D$ as in \eqref{sl(3,C)} are
\begin{equation*}
\begin{aligned}
de^1= &\, -\frac{3}{4}e^{16}+a_{1,1} e^{17}-a_{1,2} e^{27}+a_{3,1} e^{37}-a_{3,2} e^{47}+a_{5,1} e^{57}-a_{5,2} e^{67},\\
de^2= &\, \frac{3}{4}e^{26}+e^{35}+a_{1,2} e^{17}+a_{1,1} e^{27}+a_{3,2} e^{37}+a_{3,1} e^{47}+a_{5,2} e^{57}+a_{5,1} e^{67},\\
de^3= &\, \frac{1}{4}e^{36}+e^{45}+a_{1,3} e^{17}-a_{1,4} e^{27}+a_{3,3} e^{37}-a_{3,4} e^{47}+a_{5,3} e^{57}-a_{5,4} e^{67},\\
de^4= &\, -\frac{1}{4}e^{46}+e^{15}+a_{1,4} e^{17}+a_{1,3} e^{27}+a_{3,4} e^{37}+a_{3,3} e^{47}+a_{5,4} e^{57}+a_{5,3} e^{67},\\
de^5= &\, \frac{1}{2}e^{56}+a_{1,5} e^{17}-a_{1,6} e^{27}+a_{3,5} e^{37}-a_{3,6} e^{47}+\left(-a_{1,1}-a_{3,3}\right) e^{57}\\
           &\, +\left(a_{1,2}+a_{3,4}\right) e^{67},\\
de^6= &\, a_{1,6} e^{17}+a_{1,5} e^{27}+a_{3,6} e^{37}+a_{3,5} e^{47}+\left(-a_{1,2}-a_{3,4}\right) e^{57}\\
           &\, +\left(-a_{1,1}-a_{3,3}\right) e^{67},\\
de^7= &\, 0.\end{aligned}
\end{equation*}
The condition $d^2=0$ is satisfied if and only if all the parameters vanish. Thus $\mathfrak{g}$ is the direct sum
\begin{align*}
A_{6,71}^{-\frac{3}{2}}\oplus \mathbb{R}e_7 = & (-\frac{3}{4}e^{16}, \frac{3}{4}e^{26}+e^{35}, \frac{1}{4}e^{36}+e^{45}, -\frac{1}{4}e^{46}+e^{15}, \frac{1}{2}e^{56},0,0),
\end{align*}
and the $\mathrm{G}_2$ form \eqref{G2form} is closed.
\medskip


\begin{itemize}
\item \textbf{ $\frh=N_{6,13}^{0,-2,0,-2}$}
\end{itemize}
The structure equations of 
\begin{equation*}
N_{6,13}^{0,-2,0,-2}\oplus_{D} \mathbb{R}e_7,
\end{equation*}
with $D$ as in \eqref{sl(3,C)} are
\begin{equation*}
\begin{aligned}
de^1= &\, -2  \cdot 3^{-1/6}e^{16}+a_{1,1} e^{17}-a_{1,2} e^{27}+a_{3,1} e^{37}-a_{3,2} e^{47}+a_{5,1} e^{57}-a_{5,2} e^{67},\\
de^2= &\, 2 \cdot 3^{-1/6}e^{26}+e^{35}+a_{1,2} e^{17}+a_{1,1} e^{27}+a_{3,2} e^{37}+a_{3,1} e^{47}+a_{5,2} e^{57}+a_{5,1} e^{67},\\
de^3= &\, 3^{-1/6}e^{36}+ 3^{5/6}e^{45}+a_{1,3} e^{17}-a_{1,4} e^{27}+a_{3,3} e^{37}-a_{3,4} e^{47}+a_{5,3} e^{57}-a_{5,4} e^{67},\\
de^4= &\, a_{1,4} e^{17}+a_{1,3} e^{27}+a_{3,4} e^{37}+a_{3,3} e^{47}+a_{5,4} e^{57}+a_{5,3} e^{67},\\
de^5= &\, 3^{-1/6}e^{34}+ 3^{-1/6}e^{56}+a_{1,5} e^{17}-a_{1,6} e^{27}+a_{3,5} e^{37}-a_{3,6} e^{47}+\left(-a_{1,1}-a_{3,3}\right) e^{57}\\
           &\, +\left(a_{1,2}+a_{3,4}\right) e^{67},\\
de^6= &\, a_{1,6} e^{17}+a_{1,5} e^{27}+a_{3,6} e^{37}+a_{3,5} e^{47}+\left(-a_{1,2}-a_{3,4}\right) e^{57}\\
           &\, +\left(-a_{1,1}-a_{3,3}\right) e^{67},\\
de^7= &\, 0.\end{aligned}
\end{equation*}
Thus, imposing the condition $d^2=0$ we obtain that the unique derivation of $N_{6,13}^{0,-2,0,-2}$ of the form \eqref{sl(3,C)} is given by the null matrix. Hence, $\mathfrak{g}$ is the direct sum
\begin{align*}
N_{6,13}^{0, -2, 0, -2}\oplus \mathbb{R}e_7 = & (-2 \cdot 3^{-1/6}e^{16}, 2  \cdot3^{-1/6}e^{26},  3^{-1/6}e^{36}+ 3^{5/6}e^{45}, 0, \\
& \, 3^{-1/6}e^{34}+ 3^{-1/6}e^{56} ,0,0).
\end{align*}
Therefore, the $\mathrm{G}_2$ form \eqref{G2form} is closed.
\medskip

\end{proof}
\end{proposition}
\end{section}


\bigskip


\medskip

I would like to thank to Jorge Lauret and David Andriot for useful remarks and suggestions to improve the present work.
\newpage


\begin{section}{Appendix}

\begin{table}[h]
\caption{Lie algebras endowed with a closed $\mathrm{G}_2$-structure obtained in Proposition \ref{proposition}}
\begin{tabular}{cl}
\hline
$\mathfrak{g}$ & Structure equations\\\hline \hline
$\mathfrak{e}(1,1)\oplus \mathfrak{e}(1,1) \oplus_D \mathbb{R} e_7$ &$(0,0,-e^{14}+a_{3,3}e^{37},-e^{13}+a_{3,3}e^{47},e^{25}-a_{3,3}e^{57},$ \\
& $-e^{26}-a_{3,3}e^{67})$\\\hline
$(\frg_{5,1}\oplus \mathbb{R})\oplus_D \mathbb{R} e_7$  & $(0,0,a_{1,3}e^{17}+a_{3,5}e^{57},e^{15}+a_{1,3}e^{27}+a_{3,5}e^{67},a_{1,5}e^{17}+a_{3,5}e^{37}$,\\
& $e^{13}+a_{1,5}e^{27}+a_{3,5}e^{47},0)$\\\hline
$(\frg_{5,7}^{-1,-1,1}\oplus \mathbb{R})\oplus_D \mathbb{R} e_7 $ & $(-e^{15}+a_{1,3}e^{17}+a_{1,1}e^{37},e^{25}+a_{1,1}e^{17}+a_{3,1}e^{37}$,\\
& $-e^{35}+a_{1,3}e^{17}+a_{1,1}e^{37},e^{45}+a_{1,3}e^{27}-a_{1,1}e^{47},0,0,0)$\\\hline
$\big(\mathfrak{g}_{5,17}^{\alpha,-\alpha,1}\oplus\mathbb{R}\big)\oplus_{D} \mathbb{R}e_7$ & $(\alpha e^{15}+e^{35}-a_{1,3}e^{37}, -\alpha e^{25}+ e^{45}-a_{1,3}e^{47}$, \\
 & $-e^{15}+\alpha e^{35}+a_{1,3}e^{17},  -e^{25}- \alpha e^{45}+a_{1,3}e^{27}, 0, 0, 0) $ \\ \hline
$\mathfrak{g}_{6,N3}\oplus_{D} \mathbb{R}e_7 $ & $(\frac{a_{1,3}}{2}e^{37}-a_{1,5}e^{57}, e^{35}+\frac{a_{1,3}}{2} e^{47}-a_{1,5}e^{67},  a_{1,3}e^{17}+2a_{3,5}e^{57}$, \\
 & $a_{1,3}e^{27}+2a_{3,5}e^{67}, a_{1,5}e^{17}+a_{3,5}e^{37},  e^{13}+a_{1,5} e^{27}+a_{3,5}e^{47}, 0)$\\\hline
 $\mathfrak{g}_{6,38}^0\oplus \mathbb{R}e_7 $ & $(2e^{36}, 0, -e^{26}, -e^{26}+e^{25}, -e^{23}-e^{24}, e^{23},0)$\\\hline
$\frg_{6,54}^{0,-1}\oplus \mathbb{R}e_7$ & $ (e^{16}+e^{45},-e^{26},-e^{36}+e^{25},e^{46},0,0,0)$\\\hline
$\mathfrak{g}_{6,118}^{0,-1,-1}\oplus_D \mathbb{R}e_7$ & $(-e^{15}+e^{36}-a_{1,3}e^{37}, e^{46}+e^{25}-a_{1,3}e^{47}, -e^{16}-e^{35}+a_{1,3}e^{17}$,\\
& \, $e^{45}-e^{26}-a_{1,3}e^{27}, 0, 0, 0)$\\\hline
$A_{6,13}^{-\frac{2}{3},\frac{1}{3},-1}\oplus_D \mathbb{R}e_7 $ & $\big(-\frac{1}{4}e^{14}-e^{23}+a_{1,1}e^{17}, \frac{1}{4}e^{24}+a_{1,1}e^{27}, -\frac{1}{2}e^{34}, 0,$ \\ 
& \,$-\frac{3}{4} e^{45}-a_{1,1}e^{57}, \frac{3}{4} e^{46}-a_{1,1}e^{67}, 0 \big)$\\\hline
$A_{6,54}^{2,1}\oplus_D \mathbb{R}e_7 $ & $\big(-\frac{1}{2}e^{15}+a_{1,1}e^{17}+a_{3,1}e^{37}, \frac{1}{2}e^{25}+e^{16}+a_{1,1}e^{27}+a_{3,1}e^{47}, $\\
& \, $-\frac{1}{2}e^{35}+a_{1,3}e^{17}-a_{1,1}e^{37}, -\frac{1}{2}e^{45}+e^{36}+a_{1,3}e^{27}-a_{1,1}e^{47},$\\
& $0, e^{56}, 0 \big)$\\\hline
$A_{6,70}^{\alpha,\frac{\alpha}{2}}\oplus_D \mathbb{R}e_7 $ & $\big(-\frac{1}{2}e^{15}+\frac{1}{\alpha}e^{35}+e^{26}-a_{1,3}e^{37}, \frac{1}{2}e^{25}+\frac{1}{\alpha}e^{45}-a_{1,3}e^{47},$\\
& \,
$-\frac{1}{\alpha}e^{15}-\frac{1}{2}e^{35}+e^{46}+a_{1,3}e^{17}, -\frac{1}{\alpha}e^{25}+\frac{1}{2}e^{45}+a_{1,3}e^{27},$\\
& $0, e^{56}, 0 \big)
$\\\hline
$A_{6,71}^{-\frac{3}{2}}\oplus \mathbb{R}e_7 $ & $(-\frac{3}{4}e^{16}, \frac{3}{4}e^{26}+e^{35}, \frac{1}{4}e^{36}+e^{45}, -\frac{1}{4}e^{46}+e^{15}, \frac{1}{2}e^{56},0,0)$\\\hline
$N_{6,13}^{0, -2, 0, -2}\oplus \mathbb{R}e_7 $ & $(-2 \cdot 3^{-1/6}e^{16}, 2  \cdot3^{-1/6}e^{26},  3^{-1/6}e^{36}+ 3^{5/6}e^{45}, 0,$ \\
& \, $3^{-1/6}e^{34}+ 3^{-1/6}e^{56} ,0,0)$\\\hline
\end{tabular}
\end{table}

\end{section}



\bigskip

\small\noindent Universidad de Zaragoza, Departamento de Matem\'aticas,
C/Pedro Cerbuna, 12 - 50009 Zaragoza \\
\texttt{vmanero@unizar.es}



\end{document}